
\def\LaTeX{%
  \let\Begin\begin
  \let\End\end
  \let\salta\relax
  \let\finqui\relax
  \let\futuro\relax}

\def\UK{\def\our{our}\let\sz s}
\def\USA{\def\our{or}\let\sz z}

\UK



\LaTeX

\USA


\salta

\documentclass[twoside,12pt]{article}
\setlength{\textheight}{24cm}
\setlength{\textwidth}{16cm}
\setlength{\oddsidemargin}{2mm}
\setlength{\evensidemargin}{2mm}
\setlength{\topmargin}{-15mm}
\parskip2mm


\usepackage[usenames,dvipsnames]{color}
\usepackage{amsmath}
\usepackage{amsthm}
\usepackage{amssymb}
\usepackage[mathcal]{euscript}
\usepackage{cite}

%
%


\definecolor{viola}{rgb}{0.3,0,0.7}
\definecolor{ciclamino}{rgb}{0.5,0,0.5}
\definecolor{rosso}{rgb}{0.8,0,0}

\def\gianni #1{{\color{red}#1}}
\def\pier #1{{\color{rosso}#1}}
\def\juerg #1{{\color{Green}#1}}
\def\Gianni #1{{\color{blue}#1}}

\def\pier #1{#1}
\def\juerg #1{#1}
\def\gianni #1{#1}
\def\Gianni #1{#1}



\bibliographystyle{plain}


%

\finqui

\def\Beq{\Begin{equation}}
\def\Eeq{\End{equation}}
\def\Bsist{\Begin{eqnarray}}
\def\Esist{\End{eqnarray}}

\def\Bthm{\Begin{theorem}}
\def\Ethm{\End{theorem}}

\def\Brem{\Begin{remark}\rm}
\def\Erem{\End{remark}}

\def\Bcenter{\Begin{center}}
\def\Ecenter{\End{center}}
\let\non\nonumber




\def\step #1 \par{\medskip\noindent{\bf #1.}\quad}


\def\Lip{Lip\-schitz}
\def\Holder{H\"older}

\def\aand{\quad\hbox{and}\quad}

\def\lhs{left-hand side}
\def\rhs{right-hand side}
\def\sfw{straightforward}


\def\generaliz{generali\sz}

\def\regulariz{regulari\sz}


\def\multibold #1{\def\arg{#1}%
  \ifx\arg\pto \let\next\relax
  \else
  \def\next{\expandafter
    \def\csname #1#1#1\endcsname{{\bf #1}}%
    \multibold}%
  \fi \next}

\def\pto{.}

\def\multical #1{\def\arg{#1}%
  \ifx\arg\pto \let\next\relax
  \else
  \def\next{\expandafter
    \def\csname cal#1\endcsname{{\cal #1}}%
    \multical}%
  \fi \next}


\def\multimathop #1 {\def\arg{#1}%
  \ifx\arg\pto \let\next\relax
  \else
  \def\next{\expandafter
    \def\csname #1\endcsname{\mathop{\rm #1}\nolimits}%
    \multimathop}%
  \fi \next}

\multibold
qwertyuiopasdfghjklzxcvbnmQWERTYUIOPASDFGHJKLZXCVBNM.

\multical
QWERTYUIOPASDFGHJKLZXCVBNM.

\multimathop
ad dist div dom meas sign supp .


\def\accorpa #1#2{\eqref{#1}--\eqref{#2}}
\def\Accorpa #1#2 #3 {\gdef #1{\eqref{#2}--\eqref{#3}}%
  \wlog{}\wlog{\string #1 -> #2 - #3}\wlog{}}


\def\separa{\noalign{\allowbreak}}

\def\graffe #1{\mathopen\{#1\mathclose\}}

\def\<#1>{\mathopen\langle #1\mathclose\rangle}
\def\norma #1{\mathopen \| #1\mathclose \|}

\def\[#1]{\mathopen\langle\!\langle #1\mathclose\rangle\!\rangle}

\def\iot {\int_0^t}
\def\ioT {\int_0^T}
\def\intQt{\int_{Q_t}}
\def\intQ{\int_Q}
\def\iO{\int_\Omega}

\def\dt{\partial_t}
\def\dn{\partial_\nu}

\def\cpto{\,\cdot\,}

\def\checkmmode #1{\relax\ifmmode\hbox{#1}\else{#1}\fi}
\def\aeO{\checkmmode{a.e.\ in~$\Omega$}}
\def\aeQ{\checkmmode{a.e.\ in~$Q$}}

\def\aeS{\checkmmode{a.e.\ on~$\Sigma$}}

\def\aaQ{\checkmmode{for a.a.~$(x,t)\in Q$}}


\def\erre{{\mathbb{R}}}

\def\enne{{\mathbb{N}}}




\def\genspazio #1#2#3#4#5{#1^{#2}(#5,#4;#3)}
\def\spazio #1#2#3{\genspazio {#1}{#2}{#3}T0}
\def\spaziot #1#2#3{\genspazio {#1}{#2}{#3}t0}
\def\spazios #1#2#3{\genspazio {#1}{#2}{#3}s0}
\def\L {\spazio L}
\def\H {\spazio H}

\def\Lt {\spaziot L}

\def\Ls {\spazios L}
\def\C #1#2{C^{#1}([0,T];#2)}


\def\Lx #1{L^{#1}(\Omega)}
\def\Hx #1{H^{#1}(\Omega)}
\def\Wx #1{W^{#1}(\Omega)}

\def\LQ #1{L^{#1}(Q)}

\def\Luno{\Lx 1}
\def\Ldue{\Lx 2}
\def\Linfty{\Lx\infty}

\def\Huno{\Hx 1}
\def\Hdue{\Hx 2}


\def\LQ #1{L^{#1}(Q)}


\let\theta\vartheta
\let\eps\varepsilon
\let\phi\varphi

\let\TeXchi\chi                         
\newbox\chibox
\setbox0 \hbox{\mathsurround0pt $\TeXchi$}
\setbox\chibox \hbox{\raise\dp0 \box 0 }
\def\chi{\copy\chibox}



\def\coeff{1+2g(\rho)}

\def\coefft{1+2g(\rho(t))}
\def\coeffz{1+2g(\rhoz)}
\def\coeffeps{1+2g(\rhoeps)}
\def\coeffepst{1+2g(\rhoeps(t))}
\def\coeffn{1+2g(\rhon)}

\def\muz{\mu_0}
\def\rhoz{\rho_0}
\def\rhoeps{\rho_\eps}
\def\rhobar{\overline\rho}
\def\xieps{\xi_\eps}
\def\xibar{\overline\xi}
\def\mueps{\mu_\eps}

\def\mun{\mu_n}
\def\rhon{\rho_n}
\def\mubar{\overline\mu}
\def\munbar{\mubar_n}
\def\xin{\xi_n}

\def\Ldt{L^{10/3}(Q)}
\def\Lsei{\Lx6}

\def\Teps{T_\eps}

\def\normaH #1{\norma{#1}_H}

\def\betaeps{\beta_\eps}
\def\gammaeps{\gamma_\eps}
\def\tbetaeps{\widetilde\betaeps}
\def\betaz{\beta^\circ}
\def\feps{f_\eps}
\def\veps{v_\eps}

\def\Cp{C_{B,p}}
\def\CB{C_B}
\def\Cemb{C_0}
\def\Mz{M_0}
\def\calMz{\calM_0}
\def\Funo{\calF_1}
\def\Guno{\calG_1}
\def\Fdue{\calF_2}

\Begin{document}


\title{\juerg{On an application of Tikhonov's fixed point\\[0.3cm] theorem
 to a nonlocal 
  Cahn--Hilliard type\\[0.3cm] system modeling
	phase separation}}

\author{}
\date{}
\maketitle
\Bcenter
\vskip-1cm
{\large\sc Pierluigi Colli$^{(1)}$}\\
{\normalsize e-mail: {\tt pierluigi.colli@unipv.it}}\\[.25cm]
{\large\sc Gianni Gilardi$^{(1)}$}\\
{\normalsize e-mail: {\tt gianni.gilardi@unipv.it}}\\[.25cm]
{\large\sc J\"urgen Sprekels$^{(2)}$}\\
{\normalsize e-mail: {\tt sprekels@wias-berlin.de}}\\[.45cm]
$^{(1)}$
{\small Dipartimento di Matematica ``F. Casorati'', Universit\`a di Pavia}\\
{\small and Research Associate at the IMATI -- C.N.R. Pavia}
\\{\small via Ferrata 1, 27100 Pavia, Italy}\\[.2cm]
$^{(2)}$
{\small Department of Mathematics}\\
{\small Humboldt-Universit\"at zu Berlin}\\
{\small Unter den Linden 6, 10099 Berlin, Germany}\\[2mm]
{\small and}\\[2mm]
{\small Weierstrass Institute for Applied Analysis and Stochastics}\\
{\small Mohrenstrasse 39, 10117 Berlin, Germany}\\
[1cm]
{\it Dedicated to our friend Paolo Podio-Guidugli\\[.1cm]
on the occasion of his 75th birthday\\[.1cm]
with best wishes}
\Ecenter

\Begin{abstract}
\juerg{This paper investigates a nonlocal version of a model for
phase separation on an atomic lattice that was introduced by
P. Podio-Guidugli in {\it Ric.\ Mat.}\ {\bf 55}\ (2006) 105-118.
The model consists of an initial-boundary value problem for
a nonlinearly coupled system of two partial differential equations
governing the evolution of an order parameter $\rho$ and the
chemical potential $\mu$. Singular contributions to the local 
free energy in the form of logarithmic or double-obstacle potentials
are admitted. In contrast to the local model, which 
was studied by P. Podio-Guidugli and the present authors
in a series of recent publications, \pier{in the nonlocal case} 
the equation governing the 
evolution of the order parameter contains in place
of the Laplacian 
a nonlocal expression that originates from nonlocal contributions
to the free energy and accounts for possible long-range interactions
between the atoms. It is shown that just as in the local case
the model equations are well posed, where the technique of proving
existence is entirely different: it is based on an application 
of Tikhonov's fixed point theorem in a rather unusual separable and 
reflexive Banach space.} 
\vskip3mm

\juerg{Tikhonov's fixed point theorem}.
\End{abstract}

\salta

\pagestyle{myheadings}
\newcommand\testopari{\sc Colli \ --- \ Gilardi \ --- \ Sprekels}
\newcommand\testodispari{\sc Nonlocal Cahn--Hilliard system}
\markboth{\testodispari}{\testopari}

\finqui


\section{Introduction}
\label{Intro}
\setcounter{equation}{0}

\juerg{
This paper deals with a nonlocal variant of a model for phase segregation through atom rearrangement on a 
lattice proposed in~\cite{Podio}. This model (see also \cite{CGPS} for a detailed derivation), 
which is a modification of the Fried--Gurtin approach to
phase segregation processes (cf. \cite{FG}, \cite{G}), uses an order parameter $\rho$, which in many cases
represents the (normalized) density of one of the phases and attains values in the interval $[-1,1]$, and the
chemical potential $\mu$ as unknowns. It is based on a local free energy density 
of the form} 
\begin{equation}
\label{fe1}
\juerg{
\psi=\widehat\psi(\rho,\nabla\rho,\mu)=-\mu\,\rho+F(\rho)+\frac \sigma 2\,|\nabla\rho|^2,}
\end{equation}
\juerg{where $\sigma>0$ is a physical constant and $F$ is a double-well potential, and leads
to the evolutionary system
}
\Bsist
  && 2\rho \, \dt\mu
  + \mu \, \dt\rho
  - \Delta\mu = 0
  \label{oldprima}
  \\
  && -\, \juerg{\sigma}\,\Delta\rho + F'(\rho) = \mu\,. 
   \label{oldseconda}
\Esist
The above equations are assumed to hold in $Q:=\Omega\times(0,T)$,
where $\Omega$ is a \juerg{three-dimensional} domain and $T$ is some \juerg{given} final time,
and they are complemented with proper boundary and initial conditions.
\juerg{Typical examples for the double-well potential~$F$ are given by}
\Bsist
  & F_{reg}(r) := \frac14(r^2-1)^2 \,,
  \quad r \in \erre 
  \label{regpot}
  \\
  & F_{log}(r) := ((1+r)\ln (1+r)+(1-r)\ln (1-r)) - c r^2 \,,
  \quad r \in (-1,1),
  \label{logpot}
\Esist
where $c>1$ in the latter \juerg{case} so that $F_{log}$ is nonconvex. 
The potentials \eqref{regpot} and \eqref{logpot}
are usually \juerg{referred to as the {\em classical regular\/}
and the {\em logarithmic double-well\/} potential, respectively.
These potentials are smooth in their domains, where the derivative of the latter 
becomes singular at $\,\pm 1$.
However, one can even consider nondifferentiable potentials,
where an important example is given by the so-called {\em double-obstacle} potential} given~by
\Beq
  F_{2obs}(r) := I(r) - c r^2 \,,
  \quad r \in \erre,
  \label{obspot}
\Eeq
where $c>0$ is a positive constant and $I:\erre\to[0,+\infty]$ \juerg{denotes the indicator function of~$[-1,1]$, i.e.,
we have
$I(r)=0$ if $|r|\leq1$ and $I(r)=+\infty$ otherwise.
In this case, the order parameter is subjected to the unilateral constraint $|\rho|\leq1$
and \eqref{oldseconda} should be read as a differential inclusion with $F'$ representing the
subdifferential $\partial I$ of $I$.}

\juerg{
The system \eqref{oldprima}-\eqref{oldseconda} constitutes a modification of the
Cahn--Hilliard system originally introduced in \cite{CahH}
and first studied mathematically in the seminal paper \cite{EllSh} (for a large
list of references on the original Cahn--Hilliard system, see \cite{Heida}). 
It is ill-posed, in general. In fact, it was pointed out
in \cite{CGPStorino} that an associated initial-boundary value problem with zero Neumann boundary
conditions for both $\rho$ and $\mu$ may have infinitely many smooth and even nonsmooth solutions. 
Therefore, two small regularizing parameters $\varepsilon>0$ and $\delta>0$ were introduced in 
\cite{CGPS}, which led to the regularized model equations}  
\Bsist
  && \bigl( \eps + 2\rho \bigr) \, \dt\mu
  + \mu \, \dt\rho
  - \,\Delta\mu = 0
  \label{viscprima}
  \\
  && \delta\, \dt\rho - \sigma\,\Delta\rho + F'(\rho) = \mu\,. 
   \label{viscseconda}
\Esist
\juerg{The system \eqref{viscprima}--\eqref{viscseconda}, which constitutes
a modification of the so-called {\em viscous} Cahn--Hilliard system (see \cite{EllSt}),
was analyzed in the series
of papers \cite{CGPS,CGPSdc,CGPSasy,CGPSbc,CGKPS} concerning
well-posedness, regularity, optimal control and numerical approximation. Later, the
local free energy density \eqref{fe1} was generalized to the form}
\begin{equation}
\label{fe2}
\juerg{
\psi=\widehat\psi(\rho,\nabla\rho,\mu)=-\mu\,g(\rho)+F(\rho)+\frac \sigma 2\,|\nabla\rho|^2}
\end{equation}
\juerg{with a function $g$ having suitable (see below) properties. If one puts,
without loss of generality, $\varepsilon=\delta=1$, then one obtains the
more general system}
\Bsist
  && \juerg{\bigl( \coeff \bigr) \, \dt\mu
  + \mu \, g'(\rho) \, \dt\rho
  - \Delta\mu = 0}
  \label{locprima}
  \\
  && \juerg{\dt\rho - \sigma\,\Delta\rho + F'(\rho) = \mu \,g'(\rho)},
   \label{locseconda}
\Esist
\juerg{which was investigated in the papers 
\cite{CGPSgen,CGPSmagenes,CGPStorino,CGKSvd,CGKScd}.}

\juerg{In the present paper, we replace the local term $\,\frac\sigma 2\,|\nabla\rho|^2\,$
in the local free energy density by a nonlocal expression. A prototypical case is to consider
a total free energy functional of the form
\begin{equation}
\label{fe3}
{\cal F}_{tot}[\rho]=\int_\Omega \Bigl[-\mu(x)\,g(\rho(x))+F(\rho(x))\Bigr]
\,dx \,+\,{\cal Q}[\rho]
\,, 
\end{equation}
\pier{where}
$$
{\cal Q}[\rho]:=\int_\Omega\rho(x)\int_\Omega k(|y-x|)(1-\rho(y))\,dy\,dx\,. 
$$
Employing the techniques described in, e.\,g., \cite{CGPS}, we arrive with
the variational derivative} 
\begin{equation}
\juerg{B[\rho](x)=\int_\Omega k(|y-x|)\,(1-2\,\rho(y))\,dy, \quad x\in\Omega,}
\label{B1}
\end{equation} 
\juerg{of the functional ${\cal Q}$ at the following nonlocal variant of 
the system~\accorpa{locprima}{locseconda}:}
\Bsist
  && \bigl( \coeff \bigr) \, \dt\mu
  + \mu \, g'(\rho) \, \dt\rho
  - \Delta\mu = 0
  \label{Iprima}
  \\
  && \dt\rho + B[\rho] + F'(\rho) = \mu \,g'(\rho),
   \label{Iseconda}
\Esist
\juerg{which is the system that we will investigate in the following.
However, we do not restrict ourselves to operators $B$ of the exact form given
in \eqref{B1}. In fact, we consider general operators $B$ 
acting on functions defined in~$Q$ that enjoy suitable properties. \pier{Very simple} examples that
satisfy the conditions specified below are given by time convolution 
operators of the form}
\begin{equation}
\label{iop2}
\juerg{B[\rho](x,t)=\int_0^t k(t-s)\,\rho(x,s)\,ds}
\end{equation}
\juerg{and spatial convolutions of the form}
\begin{equation}
\label{iop1}
\juerg{B[\rho](x,t)=\int_\Omega k(|y-x|)\,\rho(y,t)\,dy}
\end{equation}
\juerg{provided that the respective integral kernels $k$ are smooth enough. For
instance, the three-dimensional Newtonian potential will be admissible.
However, we will not be able to include nonlocal-in-time nonlinearities of 
hysteresis type like the classical stop, play, Prandtl-Ishlinskii
or Preisach operators (for the definitions of these hysteresis operators, see,
e.\,g., \cite{BS}).}

\pier{%
Free energies of the form \eqref{fe3} have been proposed in
\cite{GL1,GL2} and rigorously justified as macroscopic limits of
microscopic phase segregation models with particle conserving
dynamics (see also \cite{CF}).
In \cite{GL1, GL2}, starting from a microscopic model, the authors derived a macroscopic equation for phase segregation phenomena that turns out to be a nonlocal 
version of the well-known Cahn--Hilliard equation.  From the
mathematical viewpoint, this nonlocal Cahn--Hilliard equation is simpler than our system \eqref{Iprima}--\eqref{Iseconda} and has received a good 
deal of attention in the last decade (see, e.g., \cite{BH1,BH2,CKRS,Gaj,GZ,H,LP}).
Most of the theoretical results are devoted to well-posedness and some are concerned with the long-time behavior of solutions. Well-posedness and regularity issues were  analyzed for an equation with degenerate mobility and logarithmic potential in \cite{GZ} (cf. also \cite{CKRS, Gaj, GG}).  This fact required to show preliminarily that a solution stays eventually strictly away from the pure phases: the so-called separation property.
For the constant mobility case and regular potentials, some existence, uniqueness and regularity results were obtained in \cite{BH1,BH2,H}. Nonsmooth potentials are considered in \cite{CKRS}. The existence of a (connected) global attractor has been proven in \cite{FrGr} for constant mobility and singular potentials. This has been done by exploiting the energy identity obtained in \cite{CFG} as a by-product of results related to a phase separation model in binary fluids. The question whether the global attractor has finite (fractal) dimension was examined in \cite{GaGr}, where the authors proved the existence of an exponential attractor.
In \cite{ABG}, an equation that is the conserved gradient flow of a nonlocal total free energy functional is considered: the functional is characterized by a Helmholtz free energy density, which can be of logarithmic type. We finally mention the paper   
\cite{RoSp}, in which a distributed optimal control problem is studied for a nonlocal convective Cahn--Hilliard equation with degenerate mobility and singular potential in three dimensions of space.}

The present paper is organized as follows.
In the next section, we \juerg{will} list our assumptions, state the problem in a precise form
and present our results. The corresponding proofs \juerg{will be} given in the last two sections. \juerg{We remark at this place that the mathematical techniques employed
here to prove existence differ significantly from those used in, e.\,g., \cite{CGPS} to
handle the local case. Indeed, while in \cite{CGPS} a retarded argument method
was utilized, we \pier{apply} Tikhonov's fixed point theorem in a rather unusual 
functional analytic framework, namely in the space $L^2(0,T;H^1(\Omega))\cap 
L^{10/3}(Q)$.}

Now, we list a number of tools and notations \juerg{employed} throughout the paper.
We repeatedly use the Young inequalities
\Bsist
  && ab\leq \delta a^2 + \frac 1 {4\delta}\,b^2
  \aand
  ab \leq \theta a^{\frac 1\theta} + (1-\theta) b^{\frac 1{1-\theta}}
  \non
  \\
  && \quad \hbox{for every $a,b\geq 0$, \ $\delta>0$, \ and \ $\theta\in(0,1),$}
  \label{young}
\Esist
as well as the H\"older and Sobolev inequalities.
In our \juerg{three-dimensional} framework, the latter read 
\Beq
  H^1 (\Omega) \subset \Lx p
  \aand
  \norma v_{p} \leq C_\Omega \Vert v \Vert_{H^1 (\Omega)}
  \quad \hbox{for every $v\in H^1 (\Omega)$ and $p\in[1,6],$}
  \label{sobolev}
\Eeq
where $C_\Omega$~depends only on~$\Omega$,
and the embedding $H^1 (\Omega)\subset\Lx p$ is compact if $p<6$.
We also recall the continuous embedding
\Beq
  \left(\pier{\L\infty {L^2 (\Omega)} \cap \L2{H^1 (\Omega)}}\right) \subset
  \left(\Ldt \cap \L{7/3}{\Lx{14/3}}\right),
  \label{embeddings}
\Eeq
which is a consequence of the Young, Sobolev and interpolation inequalities.
In particular, there holds the inequality
\Bsist
  && \norma v_{\Ldt\cap\L2{H^1 (\Omega)}}
  \leq \Cemb \max\{ \norma v_{\L\infty {L^2 (\Omega)}} , \norma{\nabla v}_{\L2{L^2 (\Omega)}} \}
  \non
  \\
  && \quad \hbox{for every $v\in\L\infty {L^2 (\Omega)}\cap\L2{H^1 (\Omega)}\,,$}
  \label{embedding}
\Esist
where $\Cemb$ depends only on~$\Omega$ and~$T$.
Finally, in order to avoid a boring notation,
we follow a general rule to denote constants.
The small-case symbol $c$ stands for different constants which depend only
on~$\Omega$, on the final time~$T$, the shape of the nonlinearities
and on the constants and the norms of
the functions involved in the assumptions of our statements.
A~small-case symbol with a subscript like $c_\delta$
indicates that the constant might depend on the parameter~$\delta$, in addition.
Hence, the meaning of $c$ and $c_\delta$ might
change from line to line and even in the same chain of equalities or inequalities.
On the contrary, we mark precise constants that we can refer~to
by using different symbols, e.g., capital letters, like in~\eqref{sobolev}.
\juerg{Also, for the sake of brevity again,} we use the same symbol $\Phi$ 
to denote different continuous functions on $[0,+\infty)$ 
with the above dependence.

\section{Statement of the problem and results}
\label{STATEMENT}
\setcounter{equation}{0}

In this section, we describe the problem under study and give an outline of our results.
As in the introduction,
$\Omega$~is the body where the evolution takes place.
We assume $\Omega\subset\erre^3$
to~be open, bounded, connected, and smooth,
and we write $|\Omega|$ for its Lebesgue measure.
Moreover, $\Gamma$ and $\dn$ stand for
the boundary of~$\Omega$ and the outward normal derivative, respectively.
Now, we specify the assumptions on the structure of our system.
We assume~that
\Bsist
  && \beta : \erre \to 2^{\erre}
  \quad \hbox{is maximal monotone with} \quad
  0 \in \beta(0)
  \label{hpbeta}
  \\[1mm]
  && \pi : \erre \to \erre
  \quad \hbox{is \Lip\ continuous}  
  \label{hppi}
  \\
  && g: \overline{D(\beta)} \to [0,+\infty)
  \quad \hbox{is $C^2$, bounded and concave, and}
  \non
  \\  
  && \quad \hbox{$g'$ is bounded and \Lip\ continuous}.
  \label{hpg}
\Esist
\gianni{In \eqref{hpg}}, $D(\beta)$ is the effective domain of~$\beta$.
For $r\in D(\beta)$, we also use the symbol $\betaz(r)$ for
the element of $\beta(r)$ having minimum modulus
(see, e.g., \cite[p.~28]{Brezis}).
\juerg{Notice that, in the notation used in the introduction, $F'=\beta+\pi$.} \pier{Moreover, let us point out that, in the case
when $D(\beta)= \erre $, our assumption \eqref{hpg} necessarily  implies 
that $g$ is a constant function, so that our system \eqref{Iprima}--\eqref{Iseconda} completely decouples; on the other hand, the significant physical case for our model (see~\cite{Podio, CGPS, CGPSgen}) corresponds to a bounded interval for $D(\beta)\  (\subseteq [-1,1] $, say) and in this framework $g$ may be rather general.}

Next, in order to list our assumptions on the nonlocal operator~$B$
and even for a future convenience, we~set
\Bsist
  && V := \Huno, \quad
  H := \Ldue
  \aand
  W := \graffe{v\in\Hdue:\ \dn v=0}
  \label{defspazi}
  \\
  && Q_t := \Omega \times (0,t)
  \quad \hbox{for $0<t\leq T$}
  \aand
  Q := Q_T\pier{.}
  \label{defQt}
\Esist
\pier{As for the \gianni{nonlocal} operator~$B$,} we assume that it maps \gianni{$\L2H=\LQ2$} into itself, is causal, and \juerg{enjoys the following properties:}
\Bsist
  && B : \L2H \to \L2H ;
  \label{hpB}
  \\[1mm]
  && B[u]|_{Q_t} = B[v]|_{Q_t}
  \quad \hbox{whenever $u|_{Q_t} = v|_{Q_t}$,\quad for every $t\in(0,T]$;}
  \label{causal}
  \\[1mm]
  \separa
  \noalign{\vskip2pt}
  && B(L^p(Q_t)) \subset L^p(Q_t)
  \aand
  \norma{B[v]}_{L^p(Q_t)} \leq \Cp \bigl( 1 + \norma v_{L^p(Q_t)} \bigr)
  \non
  \\
  && \quad \hbox{for every $v\in\LQ p$, $t\in(0,T]$, and \juerg{$p\in \left\{2, 
  \frac{10}3,6\right\}$}};
  \label{hpBbddp}
  \\[1mm]
  \separa
  \noalign{\vskip2pt}
  && \norma{B[u]-B[v]}_{L^2(Q_t)}
  \leq \CB \norma{u-v}_{L^2(Q_t)}
  \non
  \\
  && \quad \hbox{for every $u, v\in\LQ2$ and $t\in(0,T]$;}
  \label{hpBlip}
  \\[1mm]
  && B(\L2V) \subset \L2V
  \hbox{ and, for every $v\in\L2V$ and $t\in(0,T]$},\non\\
  &&\bigl| \textstyle\intQt \nabla B[v] \cdot \juerg{\nabla v}\,dx\,ds \bigr|
  \leq \CB \bigl( 1 + \textstyle\intQt \juerg{(|v|^2+ |\nabla v|^2)} \,dx\,ds\bigr)\,.
    \label{hpBbddV}
\Esist
\Accorpa\HPstruttura hpbeta hpBbddV
In the above formulas, $\Cp$~and $\CB$ are given structural constants,
and, for any Banach space~$X$, the symbol $\norma\cpto_X$ denotes its norm.
The same notation is then used also for powers of~$X$. However, \juerg{in the following 
we simply write $\norma\cpto_p$ for the standard norm in~$\Lx p$, 
for $1\leq p\leq+\infty$.}

\vspace{5mm}
\noindent\juerg{{\bf Examples.} 
\,\,It is obvious that convolution type integral operators of
the form \eqref{iop1} or \eqref{iop2} satisfy the conditions \eqref{hpB}--\eqref{hpBbddV} provided
the kernel $k$ is smooth.  However, hysteresis operators like the classical stop, play, Prandtl-Ishlinskii 
or Preisach operators are not included. The reason for this is that 
these operators carry a nonlocal memory with respect to time.
For instance, the one-dimensional stop operator ${\cal S}$ (to take the simplest of these four operators)
only enjoys (cf. \cite{BS}) the {\em nonlocal} Lipschitz property}
\begin{equation*}
\juerg{|{\cal S}[\rho_1](t)-{\cal S}[\rho_2](t)|\,\le\,2\,\max_{0\le s\le t}
\,|\rho_1(s)-\rho_2(s)|}
\end{equation*}
\juerg{for every $t\in [0,T]$ and every $\rho_1,\rho_2\in C^0([0,T])$, and it is easily
seen that the validity of the Lipschitz condition \eqref{hpBlip} cannot
be guaranteed, in general.} 

\juerg{As a further example for which the conditions can be verified, we consider the
integral operator} 
\begin{equation}
\label{iop3}
\juerg{K[\rho](x)=\int_\Omega k(|y-x|)\,\rho(y)\,dy\,, }
\end{equation}
\juerg{which acts on functions defined in $\Omega$, and its counterpart $B$ acting on functions
defined in $Q$, which is given by \eqref{iop1}. We assume that $\,k\in C^0(0, +\infty)
\,$ satisfies the condition
}
\begin{equation}
\label{yeah}
\juerg{|k(r)|\,\le\,C_1\,r^{-\alpha} \quad\,\,\forall\,r>0 \quad\,\,\mbox{with some $C_1>0$
and $\alpha<3$}. }
\end{equation} 
\juerg{Such kernels belong to the class of {\em weakly singular} kernels. Obviously, 
\eqref{causal} holds, and since $\Omega$ is a bounded 
domain, it is well known (see, e.\,g., \cite[Sect.~8.10]{Alt}) that, for any 
$p\in (1,+\infty)$ such that $\,\alpha <\frac 3 q$, 
where $\frac 1 p +\frac 1 q =1$, the linear operator 
$K$ maps $L^p(\Omega)$ continuously (even compactly) into $C^0(\overline{\Omega})$ and thus 
into $L^p(\Omega)$. It is then an easy
exercise, using H\"older's inequality, to show that for $\,\alpha < \frac 3 2$ 
the corresponding operator $B$
satisfies all of the conditions \eqref{hpB}, \eqref{hpBbddp} and \eqref{hpBlip}. 
} 

\juerg{In order to satisfy also \eqref{hpBbddV}, we need additional assumptions, for instance, that
$k$ is continuously differentiable on $\,(0,+\infty)$ with }
\begin{equation}
\label{yeahyeah}
\juerg{|k'(r)|\,\le\, C_2\,r^{-\beta}\quad\,\,\forall\,r>0
\quad\,\,\mbox{with some $C_2>0$
and $\beta<\frac 5 2$}.}
\end{equation}
\juerg{Indeed, under this assumption we have for any $v\in L^2(0,T;V)$, using the continuity of the 
embedding $V\subset L^6(\Omega)$ and the fact that $\frac{6\beta} 5<3$,}
\begin{eqnarray}
\juerg{\left|\int_{Q_t} \nabla v\cdot \nabla B[v]\right|}
&&\juerg{\le \int_{Q_t}|\nabla v|^2\,+\,c\int_{Q_t}\left |\int_\Omega
|y-x|^{-\beta}\,|v(y,s)|\,dy\right|^2 dx\,ds}\non\\
&&\juerg{\le c\int_{Q_t}|\nabla v|^2\,+\,c\int_{Q_t}\left[\int_\Omega \frac{dy} {|y-x|^{6\beta/5}}
\right]^{5/3} \|v(s)\|_6^2\,dx\,ds}\non\\
&&\juerg{\le  c\int_{Q_t} (|v|^2+|\nabla v|^2)\,.}\non 
\end{eqnarray}
\juerg{Finally, we observe that in the important case of the (long-range) three-dimensional Newtonian potential  $\,k(r)=\frac C r$, 
for which we have $\alpha=1$ and
$\beta=2$, both \eqref{yeah} and \eqref{yeahyeah} are fulfilled.}

At this point, we can describe the problem \juerg{to be investigated}.
We assume~that
\Bsist
  &&  \muz \in V 
  \aand
  \muz \geq 0 \quad \aeO 
  \label{hpmuz}
  \\
  && \rhoz \in V , \quad
  \rhoz \in D(\beta) \quad \aeO
  \aand
  \rhoz |\betaz(\rhoz)|^{7/3} \in \Luno
  \label{hprhoz}
\Esist
\Accorpa\HPdati hpmuz hprhoz
and look for a triplet $(\mu,\rho,\xi)$ satisfying
\Bsist
  && \mu \in \H1H \cap \L\infty V \cap \L2{\Wx{2,3/2}}
  \label{regmu}
  \\
  && \mu \geq 0 \quad \aeQ
  \label{mupos}
  \\
  && \rho \in \gianni{\L\infty V}
  \aand
  \dt\rho \in \Ldt
  \label{regrho}
  \\
  && \xi \in \L2H
  \label{regxi}
\Esist
\Accorpa\Regsoluz regmu regxi
and solving the initial-boundary value problem
\Bsist
  && \bigl( \coeff \bigr) \, \dt\mu
  + \mu \, g'(\rho) \, \dt\rho
  - \Delta\mu = 0
  \quad \aeQ
  \label{prima}
  \\
  && \gianni{\dt\rho + \xi + \pi(\rho) + B[\rho] = \mu \,g'(\rho)}
  \aand
  \xi \in \beta(\rho)
  \quad \aeQ
  \label{seconda}
  \\
  && \dn\mu = 0
  \quad \aeS 
  \label{neumann}
  \\
  && \mu(0) = \muz
  \aand
  \rho(0) = \rhoz \, \pier{,} 
  \label{cauchy}
\Esist
\Accorpa\Pbl prima cauchy
where $\Sigma:=\Gamma\times(0,T)$.

Here are our results.

\Bthm
\label{Existence}
With the assumptions and notations \HPstruttura\ on the structure, 
assume \HPdati\ on the initial data.
Then, problem \Pbl\ has at least \pier{one} solution satisfying \Regsoluz.
\Ethm

\Bthm
\label{Uniqueness}
Under the assumptions of Theorem~\ref{Existence},
suppose in addition \juerg{that}
\Beq
  \muz \in \Linfty
  \aand
  \rhoz \, \bigl( \betaz(\rhoz) \bigr)^5 \in \Luno.
  \label{hpreg}
\Eeq
Then the solution to problem \Pbl\ is  unique and also satisfies
\Beq
  \mu \in \LQ\infty , \quad
  \dt\rho \in \LQ6
  \aand
  \xi \in \LQ6 .
  \label{regularity}
\Eeq
\Ethm

\Brem
\label{Moregenpbl}
One can prove at least the existence of a solution
to the more general problem obtained by replacing equation \eqref{prima}~by
\Beq
  \bigl( \coeff \bigr) \, \dt\mu
  + \mu \, g'(\rho) \, \dt\rho
  - \div \bigl( \kappa(\mu)\nabla\mu \bigr) = 0,
  \label{genprima}
\Eeq
where $\kappa:[0,+\infty)\to(0,+\infty)$
is a bounded continuous function such that $1/\kappa$ is also bounded
(like the uniformly parabolic case \juerg{discussed in} 
\pier{\cite{CGPSgen, CGPStorino, CGKScd}},
while the degenerate case \pier{also treated in 
\cite{CGPSgen}} is more delicate).
Moreover, one can insert a nonnegative source term $u$
in the \rhs\ of~\eqref{genprima}.
The \pier{requirement} $u\geq0$ is needed to ensure that $\mu\geq0$,
as one can see by testing the equation by the negative part of~$\mu$
(like in the proof of \cite[Lemma 4.1]{CGPSgen}),
and \pier{a sufficient condition} that allows to \generaliz e our results is $u\in\pier{L^\infty (Q)}$.
The introduction of such a source term
\juerg{would be necessary if a distributed control problem with the control
$u$ were to be studied.
However, as the uniqueness of the solution would be needed in order to construct
the control-to-state mapping, and since a continuous dependence 
result would have} to be proved,
one should consider the situation of \cite{CGKScd}
concerning the potential and \pier{other} data
(see, in particular,  \cite[\pier{formulas}~(2.9)--(2.12)]{CGKScd}).
\Erem


\section{Existence}
\label{EXISTENCE}
\setcounter{equation}{0}

In this section, we prove Theorem~\ref{Existence}.
Our argument relies on a fixed point argument
applied to a well-defined map $\mu\mapsto\rho\mapsto\mu$
involving equations \eqref{prima} and~\eqref{seconda}, separately.
\Gianni{In our construction, we will need two different extensions
of the function $g$ to the whole of~$\erre$.
\juerg{Although we will use the same notation in both cases},
there will be no danger of confusion,
since these extensions will be used in different steps.
Furthermore, it will \juerg{become evident that the constants 
related to these extensions},
e.g., some \Lip\ constants, depend only on the corresponding constants
 related to the original map~$g$.}

\step
The functional analytic framework

In order to make it precise, we first perform a formal estimate
and construct a basic \juerg{bound}~$\Mz$. \juerg{To this end, we}
formally multiply \eqref{prima} by~$2\mu$
and observe that
\Beq
  \bigl\{ \bigl( \coeff \bigr) \, \dt\mu
  + \mu \, g'(\rho) \, \dt\rho \bigr\} \, 2\mu
  = \dt \bigl\{ \bigl( \coeff \bigr) \mu^2 \bigr\}.
  \label{performal}
\Eeq
Hence, by integrating over~$Q_t$ with $t\in(0,T)$, we have
\Beq
  \iO \bigl( \coefft \bigr) \, |\mu(t)|^2
  + \intQt |\nabla\mu|^2
  = \iO \bigl( \coeffz \bigr) \, |\muz|^2.
  \non
\Eeq
\Gianni{The function $g$ is nonnegative.
However, for \juerg{reasons that will become evident} later on,
we want to use just the inequality $g\geq-1/3$,
i.e., $1+2g\geq1/3$.
We conclude that}
\Beq
  \max \bigl\{
    \norma\mu_{\L\infty H}^2 ,
    \norma{\nabla\mu}_{\L2H}^2
  \bigr\}
  \leq \Gianni 3 (1+2\,\sup g)\, \normaH\muz^2 .
  \label{formal}  
\Eeq
Now, we owe to the embedding inequality \eqref{embedding}
and deduce that
\Beq
  \norma\mu_{\Ldt\cap\L2V}
  \leq \Mz := \Gianni {\Cemb \,(3+6\,\sup g)^{1/2}} \, \normaH\muz \,.
  \label{defraggio}
\Eeq
\juerg{Notice that} the real number $\Mz$ just defined depends only on 
$\Omega$, $T$, $g$ and~$\muz$.
At this point, we can make the first choice we need and anticipate the next one.
The \pier{used notation} should help the reader,
since $\calM$ and $\calR$ are the spaces \juerg{in which $\mu$ and $\rho$ are sought}, respectively.
We set
\Bsist
  && \calM := \Ldt \cap \juerg{L^2(0,T;V)}
  \label{defM}
  \\
  && \calMz := \graffe{v\in\calM :\ \norma v_{\calM} \leq \Mz \enskip \hbox{and} \enskip v\geq0 \ \aeQ}
  \label{defK}
  \\
  && \calR := \pier{W^{1,10/3} (0,T; L^{10/3}  (\Omega)) \cap 
  \L\infty V } .
  \label{defR}
\Esist
The next steps are devoted to the construction of the maps
$\Funo:\calMz\to\calR$ and $\Fdue:\calR\to\calMz$.
The fixed point argument will be performed on the map $\calF:=\Fdue\circ\Funo:\calMz\to\calMz$.
The definition of $\Funo$ is based on the solution 
to the Cauchy problem for \eqref{seconda}, for a given~$\mu$,~i.e.,
\Beq
  \dt\rho + \xi + \pi(\rho) + B[\rho]
  = \mu \, g'(\rho)
  \aand
  \xi \in \beta(\rho)
  \quad \aeQ , \qquad
  \rho(0) = \rhoz \,.
  \label{secondabis}
\Eeq
We have to prove a well-posedness result.

\step
The first approximating problem

\juerg{In the following, we always assume that $\mu\in\calMz$, which implies,
in particular, that $\mu\ge 0$ almost everywhere in $Q$.}
We introduce a proper approximating problem
depending on \juerg{a} positive parameter~$\eps$.
Namely, we replace $\beta$ in \eqref{secondabis} by its Yosida \regulariz ation $\betaeps$ at level~$\eps$.
We recall that $\betaeps$ is monotone and \Lip\ continuous on $\erre$
and that $|\betaeps(r)|\leq|\betaz(r)|$ for every $r\in D(\beta)$ (see, e.g., \cite[p.~28]{Brezis}).
Next, we replace $\mu$ on the \rhs\ of \eqref{secondabis} by $\Teps(\mu)$,
where the truncation map $\Teps:\erre\to\erre$ is defined~by
\Beq
  \Teps(r) := \max \{ -1/\eps , \min \{ 1/\eps,r \} \}
  \quad \hbox{for $r\in\erre$}.
  \label{defTeps}
\Eeq
Finally, we \Gianni{temporarily} extend $g$ to the whole of~$\erre$
(\juerg{still terming the extension~$g$)
in~such a way} that
\Beq
  \hbox{$g$ is a concave $C^2$ function and $g'$ is bounded and \Lip\ continuous}.
  \label{extg}
\Eeq
We stress that \juerg{we do not  require $g$ to be globally positive
so that such an extension} actually exists.
At this point, we consider the problem of finding $\rhoeps$ such that
\Beq
  \dt\rhoeps + \betaeps(\rhoeps) + \pi(\rhoeps) + B[\rhoeps]
  = \Teps(\mu) \, g'(\rhoeps)
  \quad \aeQ
  \aand
  \rhoeps(0) = \rhoz \,.
  \label{secondaeps}
\Eeq
As it is not completely obvious that such a problem has a unique solution
(due to the presence of the nonlocal operator~$B$),
we give a proof of well-posedness.
For a while, we do not stress the dependence on~$\eps$ (which is fixed)
and often avoid the subscript~$\eps$.
Clearly, the solutions $\rhoeps\in\H1H$ of~\eqref{secondaeps} are the fixed points 
(which necessarily belong to $\H1H$) of~the nonlocal operator
$\calS:\L2H\to\L2H$ defined~by
\Beq
  \calS[v](t)
  := \rhoz
  + \iot \bigl(
    \Teps(\mu) \, g'(v) - \juerg{\gammaeps(v)} - B[v]
  \bigr) (s) \, ds\,,
  \non
\Eeq
where, for brevity, we have set $\gammaeps:=\betaeps+\pi$.
In other words, for $u,v\in\L2H$, $u=\calS v$ means that
\Beq
  u \in \H1H , \quad
  \dt u = \Teps(\mu) \, g'(v) - \gammaeps(v)- B[v]
  \aand
  u(0) = \rhoz \,.
  \label{defS}
\Eeq
We claim that some \juerg{iterate} $\calS^m$ of $\calS$ is a contraction. 
\juerg{To this end, let $v_i\in\L2H$ be given} and set $u_i:=\calS[v_i]$ for $i=1,2$. 
We immediately have, for every $t\in[0,T]$,
\Bsist
  && \frac 12 \iO |u_1(t)-u_2(t)|^2
  \leq \frac 12 \intQt |u_1- u_2|^2
  \non
  \\
  && \quad {}
  + \frac 12 \intQt \bigl|
    \Teps(\mu) \, \bigl( g'(v_1) - g'(v_2) \bigr)
    - \bigl( \gammaeps(v_1) - \gammaeps(v_2) \bigr)
    - \bigl( B[v_1] - B[v_2] \bigr) 
  \bigr|^2 .
  \non
\Esist
Now, we recall that $0\le \Teps(\mu)\leq1/\eps$, that
$g'$~and $\gammaeps$ are \Lip\ continuous, and that \eqref{hpBlip} holds.
Then, by using this and applying the Gronwall lemma, we obtain \juerg{that}
\Beq
  \norma{\calS[v_1]-\calS[v_2]}_{\Lt\infty H}^2
  \leq C \norma{v_1-v_2}_{\Lt2H}^2
  \quad \hbox{for every $t\in[0,T]$,}
  \label{stimacalT}
\Eeq
where we have marked the constant by using the capital letter~$C$ for future use.
This \juerg{inequality holds for every $v_1,v_2\in\L2H$
and will be applied} to different functions.
\juerg{We now aim to show that for arbitrary $v_1,v_2\in\L2H$ and every positive integer~$m$
it holds}
\Beq
  \norma{\calS^m[v_1]-\calS^m[v_2]}_{\Lt\infty H}^2
  \leq \frac {C^m t^{m-1}} {(m-1)!} \, \norma{v_1-v_2}_{\Lt2H}^2
  \quad \hbox{for every $t\in[0,T]$} .
  \label{stimacalTm}  
\Eeq
Indeed, \eqref{stimacalTm} with $m=1$ concides with~\eqref{stimacalT}.
By assuming that $m\geq1$ and that \eqref{stimacalTm} holds,
and applying \eqref{stimacalT} to~$\calS^m[v_i]$ and \eqref{stimacalTm} to~$v_i$, 
we deduce that
\Bsist
  && \norma{\calS^{m+1}[v_1]-\calS^{m+1}[v_2]}_{\Lt\infty H}^2
  = \norma{\calS\,\calS^m[v_1]-\calS\,\calS^m[v_2]}_{\Lt\infty H}^2
  \non
  \\
  && \leq C \norma{\calS^m[v_1]-\calS^m[v_2]}_{\Lt2H}^2
  = C \iot \normaH{(\calS^m[v_1]-\calS^m[v_2])(s)}^2 \, ds
  \non
  \\
  && \leq C \iot \frac {C^m s^{m-1}} {(m-1)!} \, \norma{v_1-v_2}_{\Ls2H}^2 \, ds
  \leq \frac {C^{m+1} t^m} {m!} \, \norma{v_1-v_2}_{\Lt2H}^2 .
  \non
\Esist
Therefore, \eqref{stimacalTm} holds for every $m$,
whence $\calS^m$ is a contraction in $\L2H$ if $m$ is large enough.
This proves that the approximating problem \eqref{secondaeps} has a unique solution
$\rhoeps\in\H1H$.

\step
Construction of the first map: existence

\juerg{Next, we will derive some priori estimates and} then let $\eps$ tend to zero.
By testing the equation in \eqref{secondaeps} by~$\rhoeps$, we obtain,
for every $t\in[0,T]$,
\Beq
  \frac 12 \iO |\rhoeps(t)|^2
  + \intQt \betaeps(\rhoeps) \, \rhoeps
  = \frac 12 \iO |\rhoz|^2
  + \intQt \bigl( \Teps(\mu) g'(\rhoeps) - \pi(\rhoeps) - B[\rhoeps] \bigr) \rhoeps \,.
  \non
\Eeq
The second term on the \lhs\ is nonnegative since 
$\betaeps$ is monotone and $\betaeps(0)=0$ due to~\eqref{hpbeta}.
As for the \rhs, we owe to the definition of~$\Teps$,
the \Lip\ continuity of~$\pi$ and \eqref{hpBbddp},
and see that
\Beq
  \intQt \bigl( \Teps(\mu) g'(\rhoeps) - \pi(\rhoeps) - B[\rhoeps] \bigr) \rhoeps 
  \leq \gianni{ c \intQt \bigl( 1 + |\rhoeps|^2 + |\mu|^2 \bigr) }.
  \non
\Eeq
By applying the Gronwall lemma, we deduce that
\Beq
  \norma\rhoeps_{\L\infty H}
  \leq c \bigl( 1 + \norma\mu_{\L2H} \bigr) 
  \leq c \bigl( 1 + \norma\mu_{\calM} \bigr) .
  \label{primastimarhoeps}
\Eeq
Furthermore, as $\rhoz\in V$ and \eqref{hpBbddV} holds,
one can prove that $\rhoeps$ belongs to $\L2V$,
so that \eqref{secondaeps} can be differentiated with respect to the space variables.
Let us skip this and just derive a bound.
We take the gradient of equation \eqref{secondaeps}, 
multiply the resulting equality by $\nabla\rhoeps$ and integrate over~$Q_t$.
We obtain that
\Bsist
  && \frac 12 \iO |\nabla\rhoeps(t)|^2
  + \intQt \betaeps'(\rhoeps) |\nabla\rhoeps|^2
  = \frac 12 \iO |\nabla\rhoz|^2
  \non
  \\
  && \quad {}
  + \intQt \Bigl(
    \Teps'(\mu) g'(\rhoeps) \nabla\mu \cdot \nabla\rhoeps
    + \Teps(\mu) g''(\rhoeps) |\nabla\rhoeps|^2
    - \pi'(\rhoeps) |\nabla\rhoeps|^2
    - \nabla B[\rhoeps] \cdot \nabla\rhoeps
  \Bigr) .
  \non
\Esist
\juerg{Both} integrals on the \lhs\ are nonnegative,
and the second term in the volume integral on the \rhs\ is nonpositive since 
$\mu\geq0$ and $g''\leq0$ (cf.~\eqref{extg}).
Moreover, $0\leq\Teps'\leq1$, $g'$~and $\pi'$ are bounded and \eqref{hpBbddV} holds.
Hence, \pier{with the help of \eqref{primastimarhoeps}} we deduce that
\Beq
  \iO |\nabla\rhoeps(t)|^2
  \leq c
  + c \intQt \bigl( 1 + |\nabla\rhoeps|^2 + |\nabla\mu|^2 \bigr) .
  \non
\Eeq
Therefore, by applying the Gronwall lemma, we conclude that
\Beq
  \norma{\nabla\rhoeps}_{\L\infty H}
  \leq c \bigl( 1 + \norma{\nabla\mu}_{\L2H} \bigr)
  \leq c \bigl( 1 + \norma\mu_{\calM} \bigr) .
  \label{secondastimarhoeps}
\Eeq
Next, as $\mu\in\Ldt$, we derive an obvious bound for 
\juerg{the family $\{\Teps(\mu)\}$} in~$\Ldt$.
Moreover, \accorpa{primastimarhoeps}{secondastimarhoeps} and the embedding \eqref{embeddings} 
imply that \juerg{$\{\rhoeps\}$} is bounded in the same space,
whence the same \juerg{follows for $\{\pi(\rhoeps)\}$ and~$\{B[\rhoeps]\}$}
 (see~\eqref{hpBbddp}).
Since $g'$ is bounded, we thus have
\Beq
  \norma{\Teps(\mu) g'(\rhoeps) - \pi(\rhoeps) - B[\rhoeps]}_{\Ldt}
  \leq c \bigl( 1 + \norma\mu_{\calM} \bigr) .
  \label{perterzastimarhoeps}
\Eeq
We term $D$ the \rhs\ of \eqref{perterzastimarhoeps} 
and set $\feps:=\Teps(\mu)g'(\rhoeps)-\pi(\rhoeps)-B[\rhoeps]$,
so that \eqref{perterzastimarhoeps} becomes $\norma\feps_{\Ldt}\leq D$.

We \juerg{can} derive a similar estimate for $\{\betaeps(\rhoeps)\}$ \pier{using the following strategy}.
We set $\veps:=|\betaeps(\rhoeps)|^{7/3}\sign\betaeps(\rhoeps)$
(with $\sign 0:=0$)
and observe that $\veps\in\L\infty H$ since 
$\betaeps$ is \Lip\ continuous, 
$\rhoeps\in\L\infty V$ and $V\subset\Lx{14/3}$ by~\eqref{sobolev}.
Then, we multiply \eqref{secondaeps} by $v_\eps$
and integrate over~$Q$.
We have 
\Beq
  \iO \tbetaeps(\rho(T))
  + \intQ |\betaeps(\rho)|^{10/3}
  = \iO \tbetaeps(\rhoz)
  + \intQ \feps |\betaeps(\rhoeps)|^{7/3}\sign\betaeps(\rhoeps), 
  \non
\Eeq
where we have set
\Beq
  \tbetaeps(r) := \int_0^r |\betaeps(s)|^{7/3} \sign\betaeps(s) \, ds
  \quad \hbox{for $r\in\erre$}.
  \non
\Eeq
\pier{Note that} $\tbetaeps$ is nonnegative,
$|\tbetaeps(r)|\leq|r|\,|\betaeps(r)|^{7/3}\leq|r|\,|\betaz(r)|^{7/3}$ for every $r\in D(\beta)$
and \eqref{hprhoz} holds\pier{. Then,}
by applying the second Young inequality~\eqref{young} with $\theta=3/10$, we deduce that
\Bsist
  && \intQ |\betaeps(\rhoeps)|^{10/3}
  \leq \iO |\rhoz| \, |\betaz(\rhoz)|^{7/3}
  + \frac 3{10} \intQ |\feps|^{10/3} + \frac 7{10} \intQ |\betaeps(\rhoeps)|^{10/3}
  \non
  \\
  && \leq c + \frac 3{10} \, D^{10/3} + \frac 7{10} \intQ |\betaeps(\rhoeps)|^{10/3},
  \non
\Esist
whence immediately \juerg{$$\intQ|\betaeps(\rhoeps)|^{10/3}\leq c+D^{10/3}.$$}
We conclude that
\Beq
  \norma{\betaeps(\rhoeps)}_{\Ldt} \leq c \bigl( 1 + \norma\mu_{\calM} \bigr) \,.
  \label{terzastimarhoeps}
\Eeq
By comparison in \eqref{secondaeps} and thanks to our previous estimates, we easily infer 
that \juerg{also}
\Beq
  \norma{\dt\rhoeps}_{\Ldt} \leq c \bigl( 1 + \norma\mu_{\calM} \bigr) \,.
  \label{quartastimarhoeps}
\Eeq
At this point, it is \sfw\ to deduce that (for a subsequence) 
\Bsist
  & \rhoeps \to \rho
  & \quad \hbox{weakly star in $\L\infty V$}
  \non
  \\
  & \dt\rhoeps \to \dt\rho
  & \quad \hbox{weakly in $\Ldt$}
  \non
  \\
  & \betaeps(\rhoeps) \to \xi
  & \quad \hbox{weakly in $\Ldt$}.
  \non
\Esist
Moreover, \juerg{$\{\rhoeps\}$} converges to $\rho$ strongly in $\C0{\Lx p}$
for $p<6$, due to the compact embedding $V\subset\Lx p$ 
(see, e.g., \cite[Sect.~8, Cor.~4]{Simon}).
In particular, $\rho(0)=\rhoz$.
We also derive that \juerg{$\{B\pier{[\rhoeps]}\}$} converges to $B[\rho]$ strongly in $\L2H$ by \eqref{hpBlip},
while \juerg{$\{g'(\rhoeps)\}$ and $\{\pi(\rhoeps)\}$ converge to $g'(\rho)$ and to~$\pi(\rho)$, respectively,
strongly in $\C0{\Lx p}$ by \Lip\ continuity.}

Next, as $\mu\in\Ldt$, we see that \juerg{$\{\Teps(\mu)\}$} converges strongly to $\mu$ in~$\LQ q$ for $q<10/3$,
so that \juerg{$\{\Teps(\mu)\,g'(\rhoeps)\}$} converges to $\mu g'(\rho)$ strongly in~$\L2H$.
Finally, since \juerg{$\{\betaeps(\rhoeps)\}$} converges to $\xi$ weakly in $\LQ2$
and \juerg{$\{\rhoeps\}$} converges to $\rho$ strongly in~$\LQ2$,
\Gianni{we can apply, e.g., \cite[Lemma~2.3, p.~38]{Barbu}
to conclude that also $\rho\in D(\beta)$ and $\xi\in\beta(\rho)$ \aeQ\ 
(whence it follows that $\rho$~takes its values in the domain of the original map~$g$ (cf.~\eqref{hpg})}.
Therefore, $(\rho,\xi)$ is a solution to the Cauchy problem~\eqref{secondabis}
with the given~$\mu$.
Notice that, just by semicontinuity, 
the \pier{a priori estimates \eqref{primastimarhoeps}, \eqref{secondastimarhoeps}, \eqref{terzastimarhoeps}, \eqref{quartastimarhoeps}} are conserved in the limit,~i.e.\pier{,}
\Beq
  \norma\rho_{\calR}
  + \norma\xi_{\Ldt}
  \leq c \bigl( 1 + \norma\mu_{\calM} \bigr) 
  \quad \hbox{for every \juerg{$\mu\in\calM_0$}}\pier{,}
  \label{stimarho}
\Eeq
\pier{with obvious definition of $\norma\cdot_{\calR}$ by \eqref{defR}.}

\step
Construction of the first map: uniqueness

Let $(\rho_i,\xi_i)$, $i=1,2$, be two \juerg{solutions to the Cauchy 
problem \eqref{secondabis} for the same $\mu\in\calM_0$}.
We write the equation for both of them
and multiply the resulting equality by $\rho:=\rho_1-\rho_2$.
Then, we integrate over~$Q_t$.
We obtain
\Bsist
  && \frac 12 \iO |\rho(t)|^2
  + \intQt (\xi_1-\xi_2) \rho
  \non
  \\
  && = \intQt \mu \bigl( g'(\rho_1) - g'(\rho_2) \bigr) \rho
  - \intQt \bigl( \pi(\rho_1) - \pi(\rho_2) \bigr) \rho
  - \intQt \bigl( B[\rho_1] - B[\rho_2] \bigr) \rho .
  \non
\Esist
The second integral on the \lhs\ is nonnegative by monotonicity.
The first one on the \rhs\ is nonpositive since 
$\mu\geq0$ and $g'$ is nonincreasing by the concavity assumption \eqref{hpg} on~$g$.
By accounting for the \Lip\ continuity of~$\pi$ and~\eqref{hpBlip},
and using the Gronwall lemma,
we conclude that $\rho=0$, i.e., $\rho_1=\rho_2$.
By comparison in \eqref{secondabis}, we \juerg{see that also} $\xi_1=\xi_2$.

At this point, we can define the first map $\Funo:\calMz\to\calR$
as well as the auxiliary map $\Guno:\calMz\to\calR$ \juerg{as follows:}
\pier{\begin{align}
    &\hbox{\emph{for $\mu\in\calMz$,
    $\, \Funo(\mu)\, $ and $\, \Guno(\mu)\, $ are the \juerg{components} $\rho$ and $\xi$}}\non\\ 
    &\hbox{\emph{of the unique solution to \eqref{secondabis}.}}
  \label{defFuno}
\end{align}
}%
By the definition of $\Funo$ and $\Guno$, \eqref{stimarho} yields
\Beq
  \norma{\Funo(\mu)}_{\calR}
  + \norma{\Guno(\mu)}_{\Ldt}
  \leq c \, \bigl( 1 + \norma\mu_{\calM} \bigr)
  \leq c \, (1+\Mz) 
  \quad \hbox{for every $\mu\in\calMz$} .
  \label{stimaFuno}
\Eeq

\step
Construction of the second map: existence

Now, \Gianni{for a given $\rho\in\calR$, we would like to consider} the initial--boundary value problem
given by \eqref{prima}, \eqref{neumann} and the first \juerg{initial condition in} \eqref{cauchy}.
\Gianni{However, the terms $g(\rho)$ and $g'(\rho)$ might be meaningless 
since $g$ is not necessarily everywhere defined (cf.~\eqref{hpg}).
Hence, we suitably extend~$g$
(in~a different way with respect to the temporary~\eqref{extg},
despite of the notation we are going to use) 
to~a $C^1$ function defined in the whole real line $\erre$
by~preserving some of the properties postulated in~\eqref{hpg}.
Namely, still writing $g$ for \juerg{this new extension for the remainder} of the present section, 
we~require that
\Bsist
  && \hbox{$g$ and $g'$ are bounded and \Lip\ continuous}
  \label{extgbis}
  \\
  && g(r) \geq -1/3 , \quad \hbox{i.e.,} \quad
  1 + 2g(r) \geq 1/3 ,
  \quad \hbox{for every $r\in\erre$}.
  \label{extgter}
\Esist
Thus, the problem we consider~is
\begin{align}
  \bigl( \coeff \bigr) \, \dt\mu
  + \mu \, g'(\rho) \, \dt\rho
  - \Delta\mu = 0 \quad \pier{\aeQ ,} \non\\
  \dn\mu = 0 \quad \pier{\aeS,} \quad
  \mu(0) = \muz \,.
  \label{primabis}
\end{align}
The equation in \eqref{primabis} is linear, but its coefficients are not smooth.
Therefore}, we \regulariz e them
by introducing $\rhoeps$ as smooth as needed and satisfying
\Bsist
  & \rhoeps \to \rho
  & \quad \hbox{strongly in $\H1H$ and weakly star in $\L\infty V$}
  \qquad
  \label{hpconvrhoeps}
  \\
  & \dt\rhoeps \to \dt\rho
  & \quad \hbox{strongly in $\Ldt$} .
  \label{hpconvdtrhoeps}
\Esist
\juerg{Without loss of generality, we} can also assume that
\Bsist
  && \norma\rhoeps_{\H1H\cap\L\infty V}
  \leq 1 + \norma\rho_{\H1H\cap\L\infty V} 
  \label{daapprox}
  \\
  && \norma{\dt\rhoeps}_{\Ldt}
  \leq 1 + \norma{\dt\rho}_{\Ldt} \,.
  \label{daapproxdt}
\Esist
\juerg{The approximating problem to be considered is then}
\begin{align}
  \bigl( \coeffeps \bigr) \, \dt\mueps
  + \mueps \, g'(\rhoeps) \, \dt\rhoeps
  - \Delta\mueps = 0 \quad \pier{\aeQ ,} \non\\
  \dn\mueps = 0 \quad \pier{\aeS,}
  \quad
  \mueps(0) = \muz \,.
  \label{primaeps}
\end{align}
\juerg{It} has a unique solution $\mueps\in\H1H\cap\L\infty V\cap\L2W$
(for \juerg{the definition of} $W$, see \eqref{defspazi}),
thanks to the regularity of the coefficients and the uniform parabolicity 
ensured by~\Gianni{\eqref{extgter}}.
Moreover, the solution is nonnegative.
Indeed, by testing the equation by~$-2\mueps^-$, 
where $\mueps^-$ is the negative part of~$\mueps$,
and using the identity
\Beq
  \bigl( 
    \bigl( \coeffeps \bigr) \, \dt\mueps
    + \mu \, g'(\rhoeps) \, \dt\rhoeps
  \bigr) (-2\mueps^-)
  = \dt \bigl( \coeffeps) |\mueps^-|^2 \bigr),
  \non
\Eeq
we immediately obtain that
\Beq
  \iO \bigl( \coeffepst \bigr) |\mueps^-(t)|^2
  + \intQt |\nabla\mueps^{\pier{-}}|^2 
  = 0 \,,
  \non
\Eeq
\Gianni{whence (cf.~\eqref{extgter}) we conclude} that $\mueps^-=0$, i.e., $\mueps\geq0$.

At this point, we perform the estimate 
\Gianni{that formally led to~\eqref{formal} and was based on the inequality \eqref{extgter}.}
Since here the argument uses $\mueps$ and~$\rhoeps$,
the calculation is completely justified.
\pier{Hence, we} obtain \pier{(cf.~\eqref{defraggio})}
\Beq
  \norma\mueps_{\Ldt\cap\L2V} \leq \Mz \,.
  \label{stimaRmueps}
\Eeq
Now, we estimate some norms of $\mueps$ in terms of suitable norms of~$\rhoeps$.
The symbol $\Phi$ denotes possibly different continuous functions,
as explained at the end of Section~\pier{\ref{Intro}}.
\juerg{First, we} test \pier{the equation in} \eqref{primaeps} by~$\dt\mueps$.
By accounting for the boundedness of $g'$
and owing to the \Holder, Young and Sobolev inequalities, we have
\Bsist
  && \intQt \bigl( \coeffeps \bigr) |\dt\mueps|^2
  + \frac 12 \iO |\nabla\mueps(t)|^2
  \non
  \\
  \separa
  && \leq \frac 12 \iO |\nabla\muz|^2
  + c \intQt \mueps \, |\dt\rhoeps| \, |\dt\mueps|
  \non
  \\
  \separa
  && \leq c + c \iot \norma{\mueps(s)}_6 \, \norma{\dt\rhoeps(s)}_3 \, \norma{\dt\mueps(s)}_2 \, ds
  \non
  \\
  \separa
  && \leq c + \frac 12 \intQt |\dt\mueps|^2
  + c \iot \norma{\dt\rhoeps(s)}_3^2 \, \juerg{\norma{\mueps(s)}_V^2} \, ds .
  \non
\Esist
At this point, we recall \Gianni{\eqref{extgter} once more} and observe that \eqref{daapproxdt} implies
\Beq
  \ioT \norma{\dt\rhoeps(s)}_3^2 \, ds
  \leq c \, \norma{\dt\rhoeps}_{\Ldt}^2
  \leq c \bigl( 1 + \norma{\dt\rho}_{\Ldt}^2 \bigr) .
  \non
\Eeq
Therefore, we can apply the Gronwall lemma and conclude that (see \eqref{defR})
\Beq
  \norma\mueps_{\H1H\cap\L\infty V}
  \leq \Phi \bigl( \norma\rho_{\calR} \bigr) .
  \label{secondastimamueps}
\Eeq
Next, we estimate the first two terms of \eqref{primaeps} 
by accounting for the \Lip\ continuity of $g$ and~$g'$.
We also use the \Holder\ and Sobolev inequalities 
and owe to \accorpa{daapprox}{daapproxdt} and~\eqref{secondastimamueps}.
We easily \juerg{see that}
\Bsist
  && \norma{(\coeffeps)\,\dt\mueps}_{\L2{\Lx{3/2}}}
  \,\leq \,c \,\norma{(1+|\rhoeps|)\,\dt\mueps}_{\L2{\Lx{3/2}}}
  \non
  \\
  && \leq \,c\, \bigl( 1 + \norma\rhoeps_{\L\infty{\Lx6}} \bigr)\, \norma{\dt\mueps}_{\L2H}
  \,\leq\, \Phi \bigl( \norma\rho_{\calR} \bigr)\,,
  \non
  \\[2mm]
  \noalign{\smallskip}
  && \norma{g'(\rhoeps)\,\mueps\,\dt\rhoeps}_{\L{10/3}{\Lx{15/7}}}
  \,\leq c \, \norma{\mueps\,\dt\rhoeps}_{\L{10/3}{\Lx{15/7}}}
  \non
  \\
  && \leq c \,\norma\mueps_{\L\infty{\Lsei}} \, \norma{\dt\rhoeps}_{\L{10/3}{\Lx{10/3}}}
  \,\leq\, \Phi \bigl( \norma\rho_{\calR} \bigr) \,.
  \non
\Esist
As $10/3>2$ and $15/7>3/2$, by \pier{comparing the terms of the equation} in \eqref{primaeps}, 
we deduce a similar bound for $\Delta\mueps$ in $\L2{\Lx{3/2}}$, whence immediately
\Beq
  \norma\mueps_{\L2{\Wx{2,3/2}}}
  \leq \Phi \bigl( \norma\rho_{\calR} \bigr)
  \label{terzastimamueps}
\Eeq
by elliptic regularity.
At this point, it is \sfw\ to see that we can let $\eps$ tend to zero
to obtain a solution $\mu$ to the problem~\eqref{primabis}.
Moreover, all \juerg{of the uniform estimates shown above are preserved 
in the limit, so that we have}
\Beq
  \mu \in \calMz 
  \aand
  \norma\mu_{\H1H\cap\L\infty V\cap\L2{\Wx{2,3/2}}}
  \leq \Phi \bigl( \norma\rho_{\calR} \bigr) .
  \label{stimamu}
\Eeq

\step
Construction of the second map: uniqueness

Next, we prove that, for a given $\rho\in\calR$, the solution $\mu$ to \eqref{primabis} is unique.
We pick \pier{two} solutions $\mu_i$, $i=1,2$,
write the equation of \eqref{primabis} for both of them,
multiply the difference by $\mu:=\mu_1-\mu_2$ and integrate over~$Q_t$.
Then, \juerg{the identity \eqref{performal} holds true} for $\mu$, and we have
\Beq
  \iO \bigl( \coefft \bigr) \, |\mu(t)|^2
  + \intQt |\nabla\mu|^2
  = 0 .
  \non
\Eeq
Thus, \Gianni{by \eqref{extgter}} \pier{we conclude that} $\mu_1=\mu_2$.

At this point, we can recall \eqref{stimamu} and define $\Fdue:\calR\to\calMz$ 
\juerg{as follows:}
\Beq
  \pier{\hbox{\emph{for \gianni{$\rho\in\calR$}, 
     $\, \Fdue(\rho)\, $ is the unique solution $\mu$ to \eqref{primabis}.}}}
  \label{defFdue}
\Eeq
\juerg{We then define $\calF$ by:}
\Beq
  \calF : \calMz \to \calMz \quad
   \hbox{\pier{\emph{is given by}}} \quad
  \calF := \Fdue \circ \Funo \,.
  \label{defF}
\Eeq

\step
The fixed point argument

We want to apply Tikhonov's fixed point theorem to~$\calF$.
To this end, we observe that the Banach space $\calM$ is both reflexive and separable
and that $\calMz$ is a nonempty, bounded, \juerg{closed} and convex subset of~$\calM$.
Hence, if we endow $\calM$ with its weak topology, then
$\calMz$~is compact, and the topology induced on it by the weak topology of $\calM$
is associated to a metric. 
Therefore, in order to apply Tikhonov's theorem,
we only need to show that $\calF$ is sequentially continuous \juerg{with respect to
the weak topology of $\calM$}.
This is equivalent to \juerg{showing that, 
for every $\mubar\in\calMz$ and every sequence $\graffe{\munbar}$
of elements of $\calMz$ converging to~$\mubar$ weakly in~$\calM$,
the sequence $\graffe{\calF(\munbar)}$  
converges to $\calF(\mubar)$ weakly in~$\calM$.}

\juerg{To this end, let $\munbar,\mubar\in\calMz$ be such that 
$\munbar\to\mubar$ weakly in~$\calM$,}
and set $\rhon:=\Funo(\munbar)$, \gianni{$\xin:=\calG_1(\munbar)$,} 
and $\mun:=\calF(\munbar)=\Fdue(\rhon)$.
Thus, we~have
\begin{align}
  &\dt\rhon + \xin + \pi(\rhon) + B[\rhon]
  = \munbar \, g'(\gianni\rhon) \non\\
  &\pier{\quad\hbox{and}\quad \xin \in \beta(\rhon)
  \quad \aeQ , \quad 
  \rhon(0) = \rhoz}
  \label{secondan}
\end{align}
and we observe that the estimate \eqref{stimarho} for $\rhon$ and $\xin$ becomes
\Beq
  \norma\rhon_{\calR}
  + \norma\xin_{\Ldt}
  \leq c \bigl( 1 + \norma\munbar_{\calM} \bigr) 
  \leq c \,.
  \non
\Eeq
Therefore, we have
\Bsist
  & \rhon \to \rho
  & \quad \hbox{weakly star in $\calR$ and strongly in $\C0H$}
  \label{j1}
  \\
  & \xin \to \xi
  & \quad \hbox{weakly in $\Ldt$}
  \label{j2}
\Esist
for some $\rho$ and $\xi$ in the above spaces, at least for a subsequence
\juerg{(which is still indexed by $n\in\enne$)}.
Now, we show that $\rho=\Funo(\mubar)$ and $\xi=\Guno(\mubar)$,~i.e.,
\Beq
  \dt\rho + \xi + \pi(\rho) + B[\rho]
  = \mubar \, g'(\rho)
  \aand
  \xi \in \beta(\rho)
  \quad \aeQ , \quad
  \rho(0) = \rhoz \,.
  \label{secondater}
\Eeq
Indeed, the above strong convergence for $\{\rhon\}$ implies
both the Cauchy condition $\rho(0)=\rhoz$ and the strong convergence in $\LQ2$
of \juerg{$\{\pi(\rhon)\}$ and $\{B[\rhon]\}$} to $\pi(\rho)$ and~$B[\rho]$, respectively,
thanks to assumptions \eqref{hppi} and~\eqref{hpBlip}.
Furthermore, we also have $\xi\in\beta(\rho)$ by, e.g., \cite[Lemma~2.3, p.~38]{Barbu}.
Finally, \juerg{$\{g'(\rhon)\}$} converges to $g'(\rho)$ strongly in $\C0H$,
and \juerg{$\{\munbar\}$} converges to $\mubar$ weakly in $\L2H$, whence
\juerg{it readily follows that
$\{\munbar\, g'(\rhon)\}$} converges to $\mubar\, g'(\rho)$ weakly in~$\LQ1$.
Therefore, the pair $(\rho,\xi)$ \juerg{is the unique solution to ~\eqref{secondater},
and thus $\rho=\Funo(\mubar)$. Moreover, the limit of 
$\{\rhon\}$ is uniquely determined, from which we may conclude that all of the above convergences, in particular \eqref{j1} and \eqref{j2}, which were initially proved to be valid only for suitable
subsequences, hold in fact true for the entire sequences.}   

At this point, by setting for convenience $\mun:=\Fdue(\rhon)$,
we have
\begin{align}
  \bigl( \coeffn \bigr) \, \dt\mun
  + \mun \, g'(\rhon) \, \dt\rhon
  - \Delta\mun = 0 \quad \pier{\aeQ}, \non\\
  \dn\mun = 0 \quad
  \pier{\aeS}, \quad  \mun(0) = \muz
  \label{priman}
\end{align}
and \eqref{stimamu} for $\mun$ becomes
\Beq
  \mun \in \calMz 
  \aand
  \norma\mun_{\H1H\cap\L\infty V\cap\L2{\Wx{2,3/2}}}
  \leq \Phi \bigl( \norma\rhon_{\calR} \bigr) .
  \non
\Eeq
As $\{\rhon\}$ converges to $\rho$ weakly \gianni{star} in~$\calR$,
$\{\mun\}$~is bounded in the above norm.
Thus, we have 
\Beq
  \mun \to \mu 
  \quad \hbox{weakly star in $\H1H\cap\L\infty V\cap\L2{\Wx{2,3/2}}$}
  \label{j3}
\Eeq
for some $\mu\in\H1H\cap\L\infty V\cap\L2{\Wx{2,3/2}}$, at least for a subsequence
\juerg{(which is still indexed by $n\in\enne$)}.
We prove that $\mu=\Fdue(\rho)$,~i.e., \pier{$\mu$ solves \eqref{primabis}.}
\juerg{Indeed, since $\{\rhon\}$ converges to $\gianni\rho$ strongly in $\C0H$,
$\{g(\rhon)\}$ and $\{g'(\rhon)\}$}
converge in the same space to $g(\rho)$ and~$g'(\rho)$, respectively, 
just by \Lip\ continuity.
\Gianni{Furthermore}, \juerg{$\{\dt\mun\}$ and $\{\dt\rhon\}$} converge to $\dt\mu$ and $\dt\rho$
at least weakly in $\L2H$.
Hence, we can pass to the limit in the equation of \eqref{priman}
and deduce the \juerg{first equality in}~\pier{\eqref{primabis}}.
On the other hand, it is clear that both the boundary condition 
and the initial condition in \pier{\eqref{primabis}}
follow from the convergence of \juerg{$\{\mun\}$} to~$\mu$. \juerg{
We conclude that $\mu=\Fdue(\rho)$, that is, $\mu$ is the unique solution 
to \pier{\eqref{primabis}}. In view of the uniqueness, we may \pier{infer} that
\eqref{j3} holds true for the entire sequence.}

\juerg{Finally, we recall that $\rho=\Funo(\mubar)$. Hence, we have
proved that $\mu=\calF(\mubar)$.}
\juerg{In conclusion, Tikhonov's theorem can be applied, and $\calF$ has at 
least one fixed point
$\mu\in {\cal M}_0$.
If we consider any such $\mu$ and the corresponding pair 
$(\rho,\xi)$ given by $\rho:=\Funo(\mu)$ and $\xi:=\Guno(\mu)$,
then the estimates \eqref{stimarho} and \eqref{stimamu} are valid},
so that the triplet $(\mu,\rho,\xi)$ is a solution to problem \Pbl\
satisfying the regularity conditions \Regsoluz.


\section{Uniqueness and regularity}
\label{UNIQUENESS}
\setcounter{equation}{0}

In this section, we prove Theorem~\ref{Uniqueness}.
More precisely, we first derive \eqref{regularity} for every solution
and then show that the solution is unique.

So, we fix any solution $(\mu,\rho,\xi)$ to problem \Pbl\ satisfying \Regsoluz.
In order to prove the regularity part of the statement,
we would like to test \juerg{\eqref{seconda}} by $\xi^5$.
As no further summability of $\xi$ besides \eqref{regxi} is known,
we approximate $\rho$ and $\xi$ as follows.

\step
First auxiliary \juerg{problem}

We observe that $B[\rho]\in\LQ2$ and consider the problem
of finding $(\rhobar,\xibar)$ satisfying \accorpa{regrho}{regxi} and
\Bsist
  && \dt\rhobar + \xibar - \pi(\rhobar) - \mu g'(\rhobar) = - B[\rho]
  \aand
  \xibar \in \beta(\rhobar)
  \quad \aeQ
  \label{secondabar}
  \\
  && \rhobar(0) = \rhoz \,.
  \label{cauchybar}
\Esist
\juerg{Obviously, $(\rho,\xi)$ is a solution 
satisfying the regularity conditions \accorpa{regrho}{regxi}. We claim
that there cannot exist another such solution. To this end,
let $(\rho_i,\xi_i)$, $i=1,2$, be two solutions satisfying \accorpa{regrho}{regxi}.
We write \eqref{secondabar} for both of them,
multiply the difference by $\rho_1-\rho_2$, 
and integrate over~$Q_t$ to obtain} 
\Bsist
  && \frac 12 \iO |\rho_1(t)-\rho_2(t)|^2
  + \intQt (\xi_1-\xi_2) (\rho_1-\rho_2)
  + \intQt (-\mu) \bigl( g'(\rho_1) - g'(\rho_2) \bigr) (\rho_1-\rho_2)
  \non
  \\
  && = - \intQt \bigl( \pi(\rho_1) - \pi(\rho_2) \bigr) (\rho_1-\rho_2) \,.
  \non
\Esist
The second and third integrals on the \lhs\ are nonnegative
since $\beta$ is monotone, $\mu$~is nonnegative, 
and $g'$ is nonincreasing (see \eqref{hpg}).
Thus, by accounting for the \Lip\ continuity of $\pi$ and applying Gronwall's lemma,
we obtain $\rho_1=\rho_2$\juerg{, which proves the claim}.

\step
Second auxiliary problem

Now, we choose $\mueps\in \juerg{L^\infty (0,T;V)}\cap\LQ\infty$ with $\mueps\geq0$ such that
\Beq
  \mueps \to \mu
  \quad \hbox{strongly in \juerg{$L^\infty(0,T;V)$}}
  \label{approxmu}
\Eeq
and consider the Cauchy problem
\Beq
  \dt\rhoeps + \betaeps(\rhoeps) + \pi(\rhoeps) - \mueps g'(\rhoeps) = -B[\rho]
  \quad \aeQ
  \aand
  \rhoeps(0) = \rhoz,
  \label{secoeps}
\Eeq
where $\betaeps$ is the Yosida \regulariz ation of~$\beta$ \juerg{and where $g$ denotes
the extension of $g$ to the whole real line $\erre$
which was introduced in Section 3 and has the properties listed in \eqref{extg}}.

Since all of the nonlinearities on the \lhs\ are \Lip\ continuous
(uniformly with respect to both space and time, since $\mueps$ is bounded) and $B[\rho]\in\LQ2$, 
problem \eqref{secoeps} has a unique solution $\rhoeps\in\H1H$.
Moreover, since $B[\rho]\in\L2V$ by~\eqref{hpBbddV},
one can easily prove that $\pier{\rhoeps , \,\dt\rhoeps} \in\L2V$ and that the equations 
can be differentiated with respect to the space variables.
Thus, we can argue as we did for the proof of~\eqref{secondastimarhoeps}
(in particular, using $\mueps\geq0$ and $g''\leq0$)
and derive a bound \juerg{for the family $\{\rhoeps\}$} in~$\L\infty V$.
At this point, it is \sfw\ to show that
\juerg{$\{(\rhoeps,\betaeps(\rhoeps))\}$} converges to some $(\rhobar,\xibar)$
weakly in $\H1H\times\L2H$
(as~$\eps$ tends to zero, at least for a subsequence)
and that $(\rhobar,\xibar)$ is a solution to problem \accorpa{secondabar}{cauchybar}.
\juerg{But, as shown in the previous step, $(\rho,\xi)$ is the unique solution to this
problem}.
Therefore, we have proved that
\Beq
  \bigl( \rhoeps , \betaeps(\rhoeps) \bigr) \to (\rho,\xi)
  \quad \hbox{weakly in $\H1H\times\L2H$}
  \label{approxrho}
\Eeq
and that the convergence holds true for the whole family.

\step
Regularity

Next, we prove that $\xi\in\LQ6$ and $\dt\rho\in\LQ6$.
To this end, we consider the solution $\rhoeps$ to \eqref{secoeps}
and first show that \juerg{the family $\{\xieps:=\betaeps(\rhoeps)\}$} is bounded in~$\LQ6$.
We write the equation in \eqref{secoeps} in the form
\Beq
  \dt\rhoeps + \xieps = \feps := \mueps g'(\rhoeps) - \pi(\rhoeps) - B[\rho]
  \aand
  \xieps = \betaeps(\rhoeps).
  \label{eqsecoeps}
\Eeq
By \eqref{approxmu} and the Sobolev inequality, 
\juerg{$\{\mueps\}$ is bounded in $\L\infty\Lsei$ and thus also} in~$\LQ6$.
Since $g'$ is bounded, also \juerg{$\{\mueps g'(\rhoeps)\}$} is bounded in~$\LQ6$.
\juerg{Moreover, $\{\pi(\rhoeps)\}$ is bounded in $L^6(Q)$, since $\pi$ is Lipschitz
continuous and $\{\rho_\epsilon\}$ is known to be bounded in $L^\infty(0,T;V)$.}
Finally, as \eqref{hpBbddp} holds with $p=6$, we derive that $B[\rho]\in\LQ6$.
Thus, $\feps\in\LQ6$ and \juerg{$\{\feps\}$ is bounded in $\LQ6$.
We skip the simple} proof that $\xieps\in\LQ6$ for $\eps>0$ 
and just derive the bound we are interested in.
We multiply \eqref{eqsecoeps} by 
$\xieps^5\in\LQ{6/5}$ and integrate over~$Q$.
By noting that $\dt\rhoeps\,\xieps^5=\dt\tbetaeps(\rhoeps)$, where
\Beq
  \tbetaeps(r) := \int_0^r \bigl( \betaeps(s) \bigr)^5 \, ds
  \quad \hbox{for $r\in\erre$.}
  \non
\Eeq
we obtain
\Beq
  \iO \tbetaeps(\rho(T))
  + \intQ \xieps^6 
  = \iO \tbetaeps(\rhoz)
  + \intQ \Gianni\feps \, \xieps^5 
  \leq \iO |\rhoz| \, |\betaeps(\rhoz)|^5
  + \intQ |\Gianni\feps| \, |\xieps|^5 .
  \non
\Eeq
As $\tbetaeps$ is nonnegative by \eqref{hpbeta}, 
$|\betaeps(r)|\leq|\betaz(r)|$ for every $r\in D(\beta)$ 
(see, e.g., \cite[p.~28]{Brezis}),
and \juerg{thanks to the second condition in} \eqref{hpreg},
we can owe to the \Holder\ and Young inequalities \Gianni{in the last term}
and deduce that $\{\xieps\}$ is bounded in~$\LQ6$.
By comparison in~\eqref{eqsecoeps}, \pier{it turns out that} 
also $\{\dt\rhoeps\}$~is bounded in~$\LQ6$.
On account of~\eqref{approxrho}, we deduce that
$\xi$ and $\dt\rho$ belong to~$\LQ6$, i.e., the second and third \juerg{assertions
in \eqref{regularity} are proved.}

\pier{In order to complete the proof of \eqref{regularity}, we observe that} 
\Beq
  \dt\rho \in \L{7/3}{\Lx{14/3}} ,
  \non
\Eeq
\pier{thus} we can account for the assumption $\muz\in\Linfty$ 
to infer that $\mu\in\LQ\infty$, i.e., the \juerg{validity of the first assertion in} \eqref{regularity}, 
by repeating the argument developed in the proof of \cite[Thm.~2.3]{CGPS},
which is based on the above summability of~$\dt\rho$.
We should remark that the quoted proof is performed with $g(r)=r$;
however, only minor changes are sufficient to \juerg{arrive at} the same conclusion in the present situation
(see also the proof of the analogous \cite[Thm.~3.7]{CGPSgen} in an even more complicated case).

\step
Uniqueness

We closely follow the proof of \cite[Thm.~2.6]{CGKSvd}
and adapt the argument developed there to our situation,
also giving the details for the reader's convenience.
Indeed, \juerg{on the one hand, some of the estimates have to be changed
due to the presence of the nonlocal operator~$B$; on the other hand, it has to be clear that 
the further assumptions that were} made in \cite{CGKSvd}
in order to prove a more complicated statement are not used here.

\juerg{To begin with, we pick two solutions $(\mu_i,\rho_i,\xi_i)$, $i=1,2$,
recalling} that $\mu_i\in\LQ\infty$ by the above proof.
We write \eqref{seconda} for both of them in the form
\Beq
  \dt\rho_i + \xi_i = w_i ,
  \label{secondai}
\Eeq
\pier{where}
\Beq
 w_i := \mu_i \, g'(\rho_i) - \pi(\rho_i) - B[\rho_i].
  \label{defwi}
\Eeq
We infer that 
\Bsist
  && (\dt\rho_1 - \dt\rho_2)
  + (\xi_1 - \xi_2) 
  = w_1 - w_2
  \quad \aeQ
  \label{quattrosei}
  \\[1mm]
  && \dt |\rho_1 - \rho_2|
  + |\xi_1 - \xi_2|
  \leq |w_1 - w_2|
  \quad \aeQ \,.
  \label{quattrosette}
\Esist
The equality \eqref{quattrosei} is an obvious consequence of \gianni{\eqref{secondai}},
while \eqref{quattrosette} can be proved 
by pointwise \juerg{multiplication of} \eqref{quattrosei} by 
\gianni{$\sign(\xi_1-\xi_2)$
in the set where $\xi_1\not=\xi_2$ 
(since either $\rho_1\not=\rho_2$ and $\sign(\rho_1-\rho_2)=\sign(\xi_1-\xi_2)$
or $\dt\rho_1=\dt\rho_2$)}
and by $\sign(\rho_1-\rho_2)$ (with $\sign 0=0$)
in the set where $\xi_1=\xi_2$.
From \eqref{quattrosei} we obtain that \aaQ\ it holds (where we avoid writing the $x$ variable for brevity)
\Beq
  \iot |\dt\rho_1(s) - \dt\rho_2(s)|\juerg{\,ds}
  \leq \iot |\xi_1(s) - \xi_2(s)| \, ds
  + \iot |w_1(s) - w_2(s)| \, ds \,,
  \non
\Eeq
while \eqref{quattrosette} yields that
\Beq
  |\rho_1(t) - \rho_2(t)|
  + \iot |\xi_1(s)-\xi_2(s)| \, ds
  \leq \iot |w_1(s) - w_2(s)| \, ds .
  \non
\Eeq
By \juerg{addition}, we deduce that
\Beq
  |\rho_1(t) - \rho_2(t)|
  + \iot |\dt\rho_1(s) - \dt\rho_2(s)| \, ds
  \leq 2 \iot |w_1(s) - w_2(s)| \, ds .
  \non
\Eeq
At this point, we recall \eqref{defwi} and infer that
\Bsist
  && |\rho_1(t) - \rho_2(t)|
  + \iot |\dt\rho_1(s) - \dt\rho_2(s)| \, ds
  \leq c \iot f(s) \, ds\,,
  \quad \hbox{where}
  \non
  \\
  && f :=
  |\rho_1 - \rho_2|
  + \gianni{|\mu_1 - \mu_2|}
  + |B[\rho_1] - B[\rho_2]| \,.
  \label{quattrotre}
\Esist
Here, and in the remainder of the proof, $c$~depends also on $\norma{\mu_i}_{\LQ\infty}, \,\,i=1,2$.
We deduce that
\Beq
  |\rho_1(t) - \rho_2(t)|^2
  \leq c \, \Bigl| \iot f(s) \, ds \Bigr|^2
  \aand
  \Bigl| \iot |\dt\rho_1(s) - \dt\rho_2(s)| \, ds \Bigr|^2
  \leq c \, \Bigl| \iot f(s) \, ds \Bigr|^2 \,,
  \non
\Eeq
whence also (by integrating over~$\Omega$ and using \juerg{Schwarz's} inequality)
\Beq
  \iO |\rho_1(t) - \rho_2(t)|^2 
  \leq c \intQt |f|^2
  \aand
  \iO \Bigl| \iot |\dt\rho_1(s) - \dt\rho_2(s)| \, ds \Bigr|^2 
  \leq c \intQt |f|^2 .
  \non
\Eeq
Now, \pier{we have that} 

$$\intQt |f|^2 \leq c \intQt \bigl( |\rho_1-\rho_2|^2 + |\mu_1-\mu_2|^2 \bigr)\,,$$
by the definition of $f$ and \eqref{hpBlip}.
Therefore, we conclude that
\Bsist
  && \iO |\rho_1(t) - \rho_2(t)|^2
  \leq D \intQt \bigl( |\rho_1-\rho_2|^2 + |\mu_1-\mu_2|^2 \bigr)
  \label{quattrodiciannove}
  \\
  && \iO \Bigl| \iot |\dt\rho_1(s) - \dt\rho_2(s)| \, ds \Bigr|^2 
  \leq c \intQt \bigl( |\rho_1-\rho_2|^2 + |\mu_1-\mu_2|^2 \bigr)\,,
  \label{circaquattrodiciotto}
\Esist
where we have marked the constant in \eqref{quattrodiciannove} for future use
by using the capital letter~$D$.

\juerg{At this point, we turn our interest to} the first equation of our system.
We write it in the form\pier{%
\Beq
  \dt u_i - \Delta\mu_i = \mu_i g'(\rho_i) \dt\rho_i,
  \label{primai}
\Eeq
where $ u_i := \bigl( 1 + 2 g(\rho_i) \bigr) \mu_i,$
for $i=1,2$. Then,} take the difference and integrate with respect to time.
With the general notation \juerg{%
$$(1*v)(t):=\iot v(s)\,ds\,, \quad \pier{t\in [0,T],}$$
}%
we have
\Beq  
  (u_1-u_2) - 1*\Delta(\mu_1-\mu_2)
  = 1*\bigl( \mu_1 g'(\rho_1) \dt\rho_1 - \mu_2 g'(\rho_2) \dt\rho_2 \bigr)\,.
  \label{intdiff}
\Eeq
Then, we multiply \eqref{intdiff} by $\mu_1-\mu_2$
and integrate over~$Q_t$.
\gianni{The contribution arising from the Laplacian is nonnegative.
Now, we owe to the assumptions \eqref{hpg} on $g'$ and to Young's inequality to obtain that
\Beq
  (u_1-u_2)(\mu_1-\mu_2)
  \geq \frac 12 \, |\mu_1-\mu_2|^2 - c \, |\rho_1-\rho_2|^2 .
  \non
\Eeq
Hence, we have}
\Bsist
  && \frac 12 \intQt |\mu_1-\mu_2|^2
  \leq c \intQt |\rho_1-\rho_2|^2
  + \intQt \bigl( 1*\bigl( \mu_1 g'(\rho_1) \dt\rho_1 - \mu_2 g'(\rho_2) \dt\rho_2 \bigr) \bigr) (\mu_1-\mu_2)
  \non
  \\
  \separa
  && \leq c \intQt |\rho_1-\rho_2|^2
  + \frac 14 \intQt |\mu_1-\mu_2|^2
  \non
  \\
  && \quad {}
  + c \intQt \Bigl( \int_0^s \bigl( |\mu_1-\mu_2| + |\rho_1-\rho_2| + |\dt\rho_1 - \dt\rho_2| \bigr)(\tau) \, d\tau \Bigl)^2
  \non
  \\
  \separa
  && \leq c \intQt |\rho_1-\rho_2|^2
  + \frac 14 \intQt |\mu_1-\mu_2|^2
  \non
  \\
  && \quad {}
  + c \intQt\int_0^s |(\mu_1-\mu_2)(\tau)|^2 \, d\tau
  + c \intQt \Bigl| \int_0^s |\dt\rho_1 - \dt\rho_2|(\tau) \, d\tau \Bigl|^2
  \non
  \\
  && = c \intQt |\rho_1-\rho_2|^2
  + \frac 14 \intQt |\mu_1-\mu_2|^2
  \non
  \\
  && \quad {}
  + c \iot \Bigl( \int_{Q_s} |(\mu_1-\mu_2)|^2 \Bigr) \, ds
  + c \intQt \Bigl| \int_0^s |\dt\rho_1 - \dt\rho_2|(\tau) \, d\tau \Bigl|^2 .
  \non
\Esist
Therefore, we \juerg{find that}
\Bsist
  && \intQt |\mu_1-\mu_2|^2
  \leq c \intQt |\rho_1-\rho_2|^2
  \non
  \\
  && \quad {}
  + c \iot \Bigl( \int_{Q_s} |(\mu_1-\mu_2)|^2 \Bigr) \, ds
  + c \intQt \Bigl| \int_0^s |\dt\rho_1 - \dt\rho_2|(\tau) \, d\tau \Bigl|^2 .
  \non
\Esist
On the other hand, an integration of \eqref{circaquattrodiciotto} over $(0,t)$ yields \juerg{the estimate}
\Bsist
  && \intQt \Bigl| \int_0^s |\dt\rho_1(s) - \dt\rho_2(\tau)| \, d\tau \Bigr|^2 
  \leq c \iot \int_{Q_s} \bigl( |\rho_1-\rho_2|^2 + |\mu_1-\mu_2|^2 \bigr)
  \non
  \\
  && \leq c \intQt |\rho_1-\rho_2|^2
  + c \iot \Bigl( \int_{Q_s} |\mu_1-\mu_2|^2 \Bigr) \, ds .
  \non
\Esist
Hence, we obtain \juerg{that}
\Beq
  (D+1) \intQt |\mu_1-\mu_2|^2
  \leq c \intQt |\rho_1-\rho_2|^2
  + c \iot \Bigl( \int_{Q_s} |(\mu_1-\mu_2)|^2 \Bigr) \, ds\,, 
  \label{eureka}
\Eeq
where $D$ is the constant appearing in \eqref{quattrodiciannove}.
At this point, we \pier{take the sum of \eqref{eureka} and \eqref{quattrodiciannove}}
\juerg{to arrive at the estimate}
\Beq
  \iO |\rho_1(t) - \rho_2(t)|^2
  + \intQt |\mu_1-\mu_2|^2
  \leq c \intQt |\rho_1-\rho_2|^2
  + c \iot \Bigl( \int_{Q_s} |(\mu_1-\mu_2)|^2 \Bigr) \, ds .
  \non
\Eeq
\juerg{Applying Gronwall's lemma, we conclude that $\rho_1=\rho_2$ and $\mu_1=\mu_2$.}
Then, a comparison in \pier{\eqref{secondai}} yields $\xi_1=\xi_2$,
and the proof is complete.


\section*{Acknowledgments}
PC and GG gratefully acknowledge some financial support from the MIUR-PRIN Grant 2010A2TFX2 ``Calculus of Variations'' 
and the GNAMPA (Gruppo Nazionale per l'Analisi Matematica, la Probabilit\`a e le loro Applicazioni) 
of INdAM (Istituto Nazionale di Alta Matematica).


\vspace{3truemm}


\Begin{thebibliography}{10}

\juerg{%
\bibitem{ABG}
H. Abels, S. Bosia, M. Grasselli, Cahn--Hilliard equation with nonlocal
singular free energies, {\em Ann. Mat. Pura Appl.~\pier{(4)}} {\bf 194} (2015)
1071-1106.}

\juerg{%
\bibitem{Alt}
H. W. Alt, ``Lineare Funktionalanalysis: Eine anwendungsorientierte
Einf\"uhrung'', Springer-Verlag, Berlin--Heidelberg, 1985.
} 

\bibitem{Barbu}
V. Barbu,
``Nonlinear semigroups and differential equations in Banach spaces'',
Noordhoff, Leyden, 
1976.

\pier{%
\bibitem{BH1} P. W. Bates, J. Han, The Neumann
    boundary problem for a nonlocal Cahn-Hilliard equation, {\it J.
    Differential Equations} \textbf{212} (2005) 235-277.
}%

\pier{%
\bibitem{BH2} P. W. Bates, J. Han, The Dirichlet
    boundary problem for a nonlocal Cahn-Hilliard equation, {\it J.
    Math. Anal. Appl.} \textbf{311} (2005) 289-312.
}%

\bibitem{Brezis}
H. Brezis,
``Op\'erateurs maximaux monotones et semi-groupes de contractions
dans les espaces de Hilbert'',
North-Holland Math. Stud.
{\bf 5},
North-Holland,
Amsterdam,
1973.

\bibitem{BS}
\juerg{M. Brokate, J. Sprekels, ``Hysteresis and phase transitions'', 
Applied Mathematical Sciences {\bf 121}, Springer, New York, 1996.}

\bibitem{CahH} 
J. W. Cahn, J. E. Hilliard, 
Free energy of a nonuniform system I. Interfacial free energy, 
{\it J. Chem. Phys.\/}
{\bf 2} (1958) 258-267.

\pier{%
\bibitem{CF} C. K. Chen, P. C. Fife,  Nonlocal models of
    phase transitions in solids, {\it Adv. Math. Sci. Appl.} \textbf{10} (2000)
    821-849.
}%

\pier{%
\bibitem{CFG} 
P. Colli, S. Frigeri, M. Grasselli, Global existence of weak solutions to a nonlocal Cahn-Hilliard-Navier-Stokes system, {\it J. Math. Anal. Appl.} 386 (2012) 428-444.     
}%

\bibitem{CGKPS}
P. Colli, G. Gilardi, P. Krej\v c\'\i, P. Podio-Guidugli, J. Sprekels,
Analysis of a time discretization scheme for a nonstandard viscous Cahn-Hilliard system,
{\it ESAIM Math. Model. Numer. Anal.} {\bf 48} (2014) 1061-1087. 

\pier{%
\bibitem{CGKSvd}
P. Colli, G. Gilardi, P. Krej\v c\'\i, J. Sprekels,
A vanishing diffusion limit in a nonstandard system of phase field equations,
{\it Evol. Equ. Control Theory} {\bf 3} (2014) 257-275.
\bibitem{CGKScd}
P. Colli, G. Gilardi, P. Krej\v c\'\i, J. Sprekels,
A continuous dependence result for a nonstandard system of phase field equations,
{\it Math. Methods Appl. Sci.} {\bf 37} (2014) 1318-1324. 
}%

\bibitem{CGPS}
P. Colli, G. Gilardi, P. Podio-Guidugli, J. Sprekels,
Well-posedness and long-time behaviour for a nonstandard viscous Cahn-Hilliard system,
{\it SIAM J. Appl. Math.} {\bf 71} (2011) 1849-1870. 
 
\pier{%
\bibitem{CGPSmagenes}
P. Colli, G. Gilardi, P. Podio-Guidugli, J. Sprekels,
Global existence for a strongly coupled Cahn-Hilliard system with viscosity,
{\it Boll. Unione Mat. Ital. (9)} {\bf 5} (2012) 495-513.
\bibitem{CGPSdc}
P. Colli, G. Gilardi, P. Podio-Guidugli, J. Sprekels,
Distributed optimal control of a nonstandard system of phase field equations,
{\it Contin. Mech. Thermodyn.} {\bf 24} (2012) 437-459.
\bibitem{CGPStorino}
P. Colli, G. Gilardi, P. Podio-Guidugli, J. Sprekels,
Continuous dependence for a nonstandard Cahn-Hilliard system with nonlinear atom mobility,
{\it Rend. Sem. Mat. Univ. Pol. Torino} {\bf 70} (2012) 27-52. 
}%

\bibitem{CGPSasy}
P. Colli, G. Gilardi, P. Podio-Guidugli, J. Sprekels,
An asymptotic analysis for a nonstandard Cahn-Hilliard system with viscosity,
{\it Discrete Contin. Dyn. Syst. Ser.~S} {\bf 6} (2013) 353-368.

\bibitem{CGPSgen}
P. Colli, G. Gilardi, P. Podio-Guidugli, J. Sprekels,
Global existence and uniqueness for a singular/degenerate Cahn-Hilliard system with viscosity,
{\it J.~Differential Equations} {\bf 254} (2013) 4217-4244. 

\bibitem{CGPSbc}
P. Colli, G. Gilardi, J. Sprekels,
Analysis and optimal boundary control of a nonstandard system of phase field equations,
{\it Milan J. Math.} {\bf 80} (2012) 119-149.

\pier{%
\bibitem{CKRS} P. Colli, P. Krej\v{c}\'{i}, E. Rocca, J. Sprekels,
    Nonlinear evolution inclusions arising from phase change
    models, {\it Czechoslovak Math. J.} \textbf{57} (2007) 1067-1098.
}%

\bibitem{EllSt} 
C. M. Elliott, A. M. Stuart, 
Viscous Cahn--Hilliard equation. II. Analysis, 
{\it J. Differential Equations\/} 
{\bf 128} (1996) 387-414.

\bibitem{EllSh} 
C. M. Elliott, S. Zheng, 
On the Cahn--Hilliard equation, 
{\it Arch. Rational Mech. Anal.\/} 
{\bf 96} (1986) 339-357.

\bibitem{FG}
\juerg{E. Fried, M. E. Gurtin, Continuum theory of thermally induced phase transitions
based on an order parameter, {\em Phys. D} {\bf 68} (1993) 326-343\pier{.}}

\pier{%
\bibitem{FrGr}
S. Frigeri, M. Grasselli, Nonlocal Cahn-Hilliard-Navier-Stokes systems with singular potentials, {\it Dyn. Partial Differ. Equ.} 9 (2012) 273Ð304.
}%

\pier{%
\bibitem{Gaj} H. Gajewski, On a nonlocal model of non-isothermal phase separation,
{\it Adv. Math. Sci. Appl.} \textbf{12} (2002) 569-586.
}%

\pier{%
\bibitem{GG} H. Gajewski, J. A. Griepentrog, A descent method for the free energy of multicomponent systems, {\it Discrete Contin. Dyn. Syst.} 15 (2006) 505-528.
}%

\pier{%
\bibitem{GZ} H. Gajewski, K. Zacharias, On a nonlocal
    phase separation model, {\it J. Math. Anal. Appl.} \textbf{286} (2003)
    11-31.
}%

\pier{%
\bibitem{GaGr}
C. G. Gal, M. Grasselli,
Longtime behavior of nonlocal Cahn-Hilliard equations, 
{\it Discrete Contin. Dyn. Syst.} 34 (2014) 145-179. 
}%

\pier{%
\bibitem{GL1} G. Giacomin, J. L. Lebowitz, Phase
    segregation dynamics in particle systems with long range
    interactions. I. Macroscopic limits, {\it J. Statist. Phys.} \textbf{87}
    (1997) 37-61.
}%

\pier{%
\bibitem{GL2} G. Giacomin, J. L. Lebowitz, Phase
    segregation dynamics in particle systems with long range
    interactions. II. Phase motion, {\it SIAM J. Appl. Math.} 
    \textbf{58} (1998) 1707-1729.}%

\bibitem{G}
\juerg{M. E. Gurtin, Generalized Ginzburg--Landau and Cahn--Hilliard equations based on a
microforce balance, {\em Phys. D} {\bf 92} (1996) 178-192\pier{.}}

\pier{%
\bibitem{H} J. Han, The Cauchy problem and steady state
    solutions for a nonlocal Cahn-Hilliard equation, {\it Electron. J. Differential Equations} \textbf{113} (2004), 9 pp.
}%

\bibitem{Heida}
\juerg{M. Heida, Existence of solutions for two types of generalized versions of the
Cahn--Hilliard equation, {\em Appl. Math.} {\bf 60} (2015) 51-90.}
    
\pier{%
\bibitem{LP} S.-O. Londen, H. Petzeltov\'{a},
Convergence of solutions of a non-local phase-field system, 
{\it Discrete Contin. Dyn. Syst. Ser. S} \textbf{4} (2011) 653-670.
}%

\bibitem{Podio}
P. Podio-Guidugli, 
Models of phase segregation and diffusion of atomic species on a lattice,
{\it Ric. Mat.} {\bf 55} (2006) 105-118.

\bibitem{RoSp}
E. Rocca, J. Sprekels,
Optimal distributed control of a nonlocal convective Cahn--Hilliard equation 
by the velocity in \juerg{three dimensions, {\it SIAM J. Control Optim.} {\bf 53} (2015) 1654-1680.}

\bibitem{Simon}
J. Simon,
{Compact sets in the space $L^p(0,T; B)$},
{\it Ann. Mat. Pura Appl.~(4)\/} 
{\bf 146} (1987) 65-96.

\End{thebibliography}

\End{document}
